
  \magnification 1200


  \newcount\fontset
  \fontset=1
  \def \dualfont#1#2#3{\font#1=\ifnum\fontset=1 #2\else#3\fi}

  \dualfont\bbfive{bbm5}{cmbx5}
  \dualfont\bbseven{bbm7}{cmbx7}
  \dualfont\bbten{bbm10}{cmbx10}

  \font \eightbf = cmbx8
  \font \eighti = cmmi8 \skewchar \eighti = '177
  \font \eightit = cmti8
  \font \eightrm = cmr8
  \font \eightsl = cmsl8
  \font \eightsy = cmsy8 \skewchar \eightsy = '60
  \font \eighttt = cmtt8 \hyphenchar\eighttt = -1

  \font \sixi = cmmi6 \skewchar \sixi = '177
  \font \sixrm = cmr6
  \font \sixsy = cmsy6 \skewchar \sixsy = '60
  \font \tensc = cmcsc10
  
  \font \titlefont = cmbx12
  \scriptfont \bffam = \bbseven
  \scriptscriptfont \bffam = \bbfive
  \textfont \bffam = \bbten

  \newskip \ttglue

  \def \eightpoint {\def \rm {\fam0 \eightrm }%
  \textfont0 = \eightrm
  \scriptfont0 = \sixrm \scriptscriptfont0 = \fiverm
  \textfont1 = \eighti
  \scriptfont1 = \sixi \scriptscriptfont1 = \fivei
  \textfont2 = \eightsy
  \scriptfont2 = \sixsy \scriptscriptfont2 = \fivesy
  \textfont3 = \tenex
  \scriptfont3 = \tenex \scriptscriptfont3 = \tenex
  \def \it {\fam \itfam \eightit }%
  \textfont \itfam = \eightit
  \def \sl {\fam \slfam \eightsl }%
  \textfont \slfam = \eightsl
  \def \bf {\fam \bffam \eightbf }%
  \textfont \bffam = \bbseven
  \scriptfont \bffam = \bbfive
  \scriptscriptfont \bffam = \bbfive
  \def \tt {\fam \ttfam \eighttt }%
  \textfont \ttfam = \eighttt
  \tt \ttglue = .5em plus.25em minus.15em
  \normalbaselineskip = 9pt
  \def \MF {{\manual opqr}\-{\manual stuq}}%
  \let \sc = \sixrm
  \let \big = \eightbig
  \setbox \strutbox = \hbox {\vrule height7pt depth2pt width0pt}%
  \normalbaselines \rm }



  \newcount \secno \secno = 0
  \newcount \stno \stno = 0
  \newcount \eqcntr \eqcntr= 0

  \def \ifn #1{\expandafter \ifx \csname #1\endcsname \relax }

  \def \track #1#2#3{\ifn{#1}\else {\tt\ [#2 \string #3] }\fi}

  \def \advseqnumbering {\global \advance \stno by 1 \global \eqcntr=0}

  \def \current {\number \secno \ifnum \number \stno = 0 \else
    .\number \stno \fi }

  \def \laberr#1#2{\message{*** RELABEL CHECKED FALSE for #1 ***}
      RELABEL CHECKED FALSE FOR #1, EXITING.
      \end}

  \def \syslabel#1#2{%
    \ifn {#1}%
      \global \expandafter 
      \edef \csname #1\endcsname {#2}%
    \else
      \edef\aux{\expandafter\csname #1\endcsname}%
      \edef\bux{#2}%
      \ifx \aux \bux \else \laberr{#1=(\aux)=(\bux)} \fi
      \fi
    \track{showlabel}{*}{#1}}

  \def \subeqmark #1 {\global \advance\eqcntr by 1
    \edef\aux{\current.\number\eqcntr}
    \eqno {(\aux)}
    \syslabel{#1}{\aux}}

  \def \eqmark #1 {\advseqnumbering
    \eqno {(\current)}\syslabel{#1}{\current}}

  \def \fcite#1#2{\syslabel{#1}{#2}\lcite{#2}}

  \def \label #1 {\syslabel{#1}{\current}}

  \def \lcite #1{(#1\track{showcit}{$\bullet$}{#1})}

  \def \cite #1{[{\bf #1}\track{showref}{\#}{#1}]}

  \def \scite #1#2{{\rm [\bf #1\track{showref}{\#}{#1}{\rm \hskip 0.7pt:\hskip 2pt #2}\rm]}}


 \def \Headlines #1#2{\nopagenumbers
    \advance \voffset by 2\baselineskip
    \advance \vsize by -\voffset
    \headline {\ifnum \pageno = 1 \hfil
    \else \ifodd \pageno \tensc \hfil \lcase {#1} \hfil \folio
    \else \tensc \folio \hfil \lcase {#2} \hfil
    \fi \fi }}

  \def \Title{\centerline{\titlefont \ucase{\titletextOne}}
    \edef\aux{\titletextTwo}\edef\bux{\null}%
    \ifx\aux\bux \else 
      \smallskip
      \centerline{\titlefont \ucase{\titletextTwo}}
      \fi}

  \def \Date #1 {\footnote {}{\eightit Date: #1.}}


  \def \ucase #1{\edef \auxvar {\uppercase {#1}}\auxvar }
  \def \lcase #1{\edef \auxvar {\lowercase {#1}}\auxvar }

  \def \section #1{\global\def \SectionName{#1}\stno = 0 \global
\advance \secno by 1 \bigskip \bigskip \goodbreak \noindent {\bf
\number \secno .\enspace #1.}\medskip \noindent \ignorespaces}

  \long \def \sysstate #1#2#3{%
    \advseqnumbering
    \medbreak \noindent 
    {\bf \current.\enspace #1.\enspace }{#2#3\vskip 0pt}\medbreak }
  \def \state #1 #2\par {\sysstate {#1}{\sl }{#2}}
  \def \definition #1\par {\sysstate {Definition}{\rm }{#1}}
  \def \remark #1\par {\sysstate {Remark}{\rm }{#1}}


  \def \proof {\medbreak \noindent {\it Proof.\enspace }}
  \def \proofend {\ifmmode \eqno \square \else \hfill \square
\looseness = -1 \medbreak \fi }

  \def \$#1{#1 $$$$ #1}
  \def \=#1{\buildrel \hbox{\sixrm #1} \over =}

  \def \Item #1{\smallskip \item {{\rm #1}}}
  \newcount \zitemno \zitemno = 0

  \def \izitem {\zitemno = 0}
  \def \zitemplus {\global \advance \zitemno by 1 \relax}
  \def \rzitem{\romannumeral \zitemno}
  \def \rzitemplus {\zitemplus \rzitem}
  \def \zitem {\Item {{\rm(\rzitemplus)}}}
  \def \zitemmark #1 {\syslabel{#1}{\rzitem}}

  \newcount \nitemno \nitemno = 0
  
  \def \nitem {\global \advance \nitemno by 1 \Item {{\rm(\number\nitemno)}}}

  \newcount \aitemno \aitemno = -1
  \def \boxlet#1{\hbox to 6.5pt{\hfill #1\hfill}}
  \def \iaitem {\aitemno = -1}
  \def \aitemconv{\ifcase \aitemno a\or b\or c\or d\or e\or f\or g\or
h\or i\or j\or k\or l\or m\or n\or o\or p\or q\or r\or s\or t\or u\or
v\or w\or x\or y\or z\else zzz\fi}
  \def \aitem {\global \advance \aitemno by 1\Item {(\boxlet \aitemconv)}}
  \def \aitemmark #1 {\syslabel{#1}{\aitemconv}}

  \newcount \footno \footno = 1
  \newcount \halffootno \footno = 1
  \def \footcntr {\global \advance \footno by 1
  \halffootno =\footno
  \divide \halffootno by 2
  $^{\number\halffootno}$}
  \def \fn#1{\footnote{\footcntr}{\eightpoint#1\par}}

  \begingroup
  \catcode `\@=11
  \global\def\eqmatrix#1{\null\,\vcenter{\normalbaselines\m@th%
      \ialign{\hfil$##$\hfil&&\kern 5pt \hfil$##$\hfil\crcr%
	\mathstrut\crcr\noalign{\kern-\baselineskip}%
	#1\crcr\mathstrut\crcr\noalign{\kern-\baselineskip}}}\,}
  \endgroup


  \font\mf=cmex10
  \def\union {\mathop{\raise 9pt \hbox{\mf S}}\limits}
  \def\inters{\mathop{\raise 9pt \hbox{\mf T}}\limits}

  \def \N {{\bf N}}
  \def \C {{\bf C}}
  \def \<{\left \langle \vrule width 0pt depth 0pt height 8pt }
  \def \>{\right \rangle }  
  
  \def \ds{\displaystyle}
  \def \and {\hbox {,\quad and \quad }}
  
  \def \labelarrow#1{\ {\buildrel #1 \over \longrightarrow}\ }
  \def \imply {\mathrel{\Rightarrow}}
  \def \for #1{,\quad \forall\,#1}
  \def \square {\hbox {$\sqcap \!\!\!\!\sqcup $}}
  
  \def \stress #1{{\it #1}\/}
  \def \inv {^{-1}}
  \def \*{\otimes}

  \newcount \bibno \bibno = 0
  \def \newbib #1{\global\advance\bibno by 1 \edef #1{\number\bibno}}

  \def \bibitem #1#2#3#4{\smallskip \item {[#1]} #2, ``#3'', #4.}

  \def \references {
    \begingroup
    \bigskip \bigskip \goodbreak
    \eightpoint
    \centerline {\tensc References}
    \nobreak \medskip \frenchspacing }

\def \caldef #1{\global \expandafter \edef \csname #1\endcsname {{\cal #1}}}

  \input pictex

  \newcount\ax
  \newcount\ay
  \newcount\bx
  \newcount\by
  \newcount\dx
  \newcount\dy
  \newcount\vecNorm
  \newcount\pouquinho \pouquinho = 200

  \def\myarrow#1#2#3#4{\arrow <0.15cm> [0.25,0.75] from #1 #2 to #3 #4 }%
  \def\morph#1#2#3#4{%
    \ax = #1
    \ay = #2
    \bx = #3
    \by = #4
    \dx = \bx \advance \dx by -\ax
    \dy = \by \advance \dy by -\ay
    \vecNorm = \dx 
    \ifnum\vecNorm<0 \vecNorm=-\vecNorm \fi
    \advance \vecNorm by \ifnum\dy>0 \dy \else -\dy \fi
    \multiply \dx by \pouquinho \divide \dx by \vecNorm
    \multiply \dy by \pouquinho \divide \dy by \vecNorm
    \advance \ax by \dx
    \advance \bx by -\dx
    \advance \ay by \dy
    \advance \by by -\dy
    \myarrow{\number\ax}{\number\ay}{\number\bx}{\number\by}}

  \def\starmap#1{{\rm star}_{#1}}
  \def\star#1#2{\starmap{#1}(#2)}

  \def\mult{\mu}
  \def\incl{j}
  \def\G{{\cal G}}
  \def\A{{\cal A}}
  \def\NA{{\cal N}}
  \def\LA{{\cal L}(\A)}
  \def\phi{\varphi}
  \def\tphi{\tilde\phi}
  
  \def\phie{\varphi_e}
  \def\tphie{\tilde\phie}
  \global\def\app{w}
  \def\d{\delta}
  \def\piplus{\pi^{\scriptscriptstyle +}}
  \def\Piplus{\Pi^{\scriptscriptstyle +}}
  \def\pizero{\pi^{\scriptscriptstyle 0}}
  \def\Pizero{\Pi^{\scriptscriptstyle 0}}
  \def\piu{\pi^u}
  \def\reduced{{\hbox{\sixrm red}}}
  \def\piur{\pi^\reduced}
  \def\CstarRed{C^*_\reduced}
  \def\iotaaRed{\iotaa^\reduced}
  \def\span{{\rm span}}
  \def\net#1,#2{\{#1\}_{#2}}
  \caldef F
  \caldef E
  \def\dcup{\mathop{\dot {\hbox{$\bigcup$}}}}
  \def\half{^{1/2}}
  
  \def\Fell{F}
  \def\tnorm#1{\def\tri{\hbox{$|\kern-1.5pt|\kern-1.5pt|$}}\tri #1 \tri}
  \def\tnormFunction{\tnorm{\kern-1.5pt\cdot\kern-1.5pt}}

  \def\Cx{P}
  \def\RiLA{\Upsilon_{\tphi}}
  \def\Bu{\widetilde B}
  \def\closure#1{\overline {#1}}
  \def\T{\tau}
  \def\vc{virtual commutant}
  \def\supp{{\rm supp}}
  \def\Ker{{\rm Ker}}
  \def\zerosp{\{0\}} 
  \def\0{\{0\}} 
  \def\imply{\ \Rightarrow \ }
  \def\cltsi{closed two-sided ideal}
  \def\clf{continuous linear functional}
  \def\gc{generalized Cartan}
  \def\Gc{Generalized Cartan}
  \def\so#1{S(#1)}
  \def\rg#1{R(#1)}
  \font\gothic=eufm10
  \font\smallgothic=eufm10 scaled 833
  \def\S{\hbox{\gothic S}}
  \def\sS{\hbox{\smallgothic S}}
  \def\SAB{\S_{A, B}}
  \def\sSAB{\sS_{A, B}}
  
  \def\prd{\null}  
  \def\clsum{\mathop{\overline {\sum}}}
  \font\smsc=cmcsc10 scaled 833
  \def\Ab{{\smsc Ab}}
  \def\Max{{\smsc Max}}
  \def\MaxPrime{\Max$'$}

  \def \newbib #1#2{\global\advance\bibno by 1 \edef #1{\number\bibno}}

  \newbib\Blackadar   {B}
  \newbib\tpa         {E1}
  \newbib\actions     {E2}
  \newbib\FMOne       {F1}
  \newbib\FMTwo       {F2}
  \newbib\CH          {H}
  \newbib\Jensen      {JT}
  \newbib\Jones       {J}
  \newbib\Kasparov    {Ka}
  \newbib\Kumjian     {Ku}
  \newbib\Lawson      {L}
  \newbib\Pedersen    {P}
  \newbib\Renault     {R}
  \newbib\TakeOne     {T1}
  \newbib\TakeTwo     {T2}
  \newbib\Watatani    {W}
  \newbib\Zettl       {Z}

  \def\titletextOne{NONCOMMUTATIVE CARTAN SUB-ALGEBRAS}
  \def\titletextTwo{OF C*-ALGEBRAS}

  \Headlines {\titletextOne \ \titletextTwo} {R.~Exel}

  \null\vskip -1cm
  \centerline{\bf \titletextOne} \smallskip
  \centerline{\bf \titletextTwo}
  \footnote{\null}
  {\eightrm 2000 \eightsl Mathematics Subject Classification:
  \eightrm 
  Primary   46L45, 
  secondary 46L55, 
            20M18. 
  }

  \centerline{\tensc 
    R.~Exel\footnote{*}{\eightpoint Partially supported by
CNPq.}}\footnote{\null}
  {\eightrm Keywords: C*-algebras, Cartan subalgebras,
  inverse semigroups, Fell bundles.}
  \Date{25 Jun 2008}

  \midinsert 
  \narrower \narrower
  \eightpoint \noindent J.~Renault has recently  found a
generalization of the caracterization of C*-diagonals obtained by
A.~Kumjian in the eighties, which in turn is a C*-algebraic version of
J.~Feldman and C.~Moore's well known Theorem on Cartan subalgebras of
von Neumann algebras.  Here we propose to give a version of Renault's
result in which the {Cartan subalgebra is not necessarily commutative}
[sic].  Instead of describing a Cartan pair as a twisted groupoid
C*-algebra we use N.~Sieben's notion of Fell bundles over inverse
semigroups which we believe should be thought of as
\stress{twisted \'etale groupoids with noncommutative unit space}.  
En passant we prove a theorem on uniqueness of
conditional expectations.
  \endinsert

\section{Introduction}
  Building on work by Feldman, Moore \cite{\FMOne,\FMTwo} and Kumjian
\cite{\Kumjian}, Renault has recently introduced a natural notion of
Cartan subalgebras of C*-algebras \cite{\Renault}.  If $B$ is a Cartan
subalgebra of the C*-algebra $A$, Renault proved the existence of an
essentially principal \'etale groupoid $G(B)$ and a 2-cocycle
$\Sigma(B)$, respectively called the \stress{Weyl groupoid} and the
\stress{Weyl twist}, such that $A$ is isomorphic to the reduced
twisted groupoid C*-algebra $C_r^*\big(G(B),\Sigma(B)\big)$, in such a
way that $B$ is carried onto the algebra of continuous functions
vanishing at infinity on the unit space of $G(B)$.  In other words,
this gives a caracterization of Cartan subalgebras in the context of
C*-algebras paralleling Feldman and Moore's caracterization of Cartan
subalgebras of von Neumann algebras.

Among the important consequences of Renault's work is the fact that
the conditional expectation from a C*-algebra to a Cartan subalgebra
is unique \scite{\Renault}{5.7}.
  Uniqueness of conditional expectations is a common and useful
phenomenon in the theory of von Neumann algebras and it holds whenever
the von Neumann subalgebra $B$ of the von Neumann algebra $A$
satisfies $B'\cap A\subseteq B$.  See \scite{\TakeTwo}{IX.4.3}. 

Renault's result on uniqueness of conditional expectations may thus be
considered as a version of this well known result,  but
the requirement that $B$ be a maximal abelian subalgebra (one of the
conditions for $B$ to be a Cartan subalgebra) is way stronger than
$B'\cap A\subseteq B$.
  Incidentally notice that $B$ is a maximal abelian subalgebra of $A$ if and
only if the following two conditions hold:
  \medskip
  \itemitem{(\Max)} $B'\cap A\subseteq B$,
  and
  \medskip
  \itemitem{(\Ab)} $B$ is abelian.

\medskip \noindent Thus, if one is hoping for a direct generalization
of the von Neumann uniqueness of conditional expectations mentioned
above to the context of C*-algebras, Renault's result should be
strengthened by removing condition (\Ab) above from the hypotheses.  

  \def\K{{\cal K}}%
  Unfortunately this is impossible.  If $A=C([0,1])\*\K$, where $\K$
is the algebra of compact operators on an infinite dimensional Hilbert
space, and $B = 1\*\K$, then the pair $(A,B)$ satisfies Renault's
axioms of Cartan pairs except that, in place of maximal abeliannes,
only axiom (\Max) holds.
  Nevertheless conditional expectations abound as each probability
measure on $[0, 1]$ gives, by means of integration, a different
conditional expectation from $C([0, 1])\*\K$ to $1\*\K$.

However the validity of (\Max) in this example sounds fishy: $B'\cap
A=\zerosp$, so it is contained in $B$ alright but, not far outside $A$
one does find \stress{\vc s} of $B$ not belonging to $B$, such as any
element of the form $f\*1$.  Technically speaking, what we mean by a
virtual commutant (Definition \fcite{VirtuComm}{9.2}) for an inclusion
``$B\subseteq A$" of C*-algebras is a $B$-bimodule map
  $$
  \phi: J\to A,
  $$
  where $J$ is a {\cltsi} of $B$.
  If $a\in B'\cap A$ and $J$ is any ideal in $B$, then $\phi(x) = ax$,
gives an example of a {\vc} defined on $J$.  Hence these should be
thought of as \stress{generalized elements} of $B'\cap A$.

However not every {\vc} is of this form.  In the above example in
which $B=1\*\K$ and $A=C([0,1])\*\K$, let $f$ be
a non-constant function and put
  $$
  \phi(1\*x) = f\*x \for x\in \K.
  $$
  Then $\phi$ is a {\vc} defined on $J = B$, but it is not of the
elementary form above.
One may then replace (\Max) above by the following stronger axiom:
  \medskip
  \itemitem{(\MaxPrime)} the range of every {\vc} is contained in $B$,

\medskip\noindent in which case our badly behaved example 
will be knocked out.
  This replacement is not at all a drastic departure from tradition
since (\MaxPrime+\Ab) is equivalent to (\Max+\Ab), and hence to maximal
abeliannes.  See \fcite{GeneralizeMaxAbel}{9.8} for a proof of this
statement.
  
We thus propose to reformulate Renault's definition of Cartan
subalgebras \scite{\Renault}{Definition 5.1} by replacing maximal
abeliannes with (\MaxPrime) above,  without requiring (\Ab).  

  I am well aware that Cartan subalgebras, wherever they have been
considered, have always been assumed to be abelian algebras.  So our
use of this term to refer to noncommutative subalgebras might sound
slightly heretic.  But since we shall provide generalizations of
well known classical results for Cartan subalgebras, we hope our
heresy will be forgiven.

Assuming that $(A, B)$ is such a {\gc} pair, and that $A$
is separable, we prove in Theorem \fcite{UniqueCondexp}{12.3} that the
conditional expectation onto $B$ is unique, thus extending Renault's
result \scite{\Renault}{5.7}.

Our second main contribution is to generalize Renault's version of
Feldman and Moore's Theorem \scite{\Renault}{5.6} but, since 
$B$ may be noncommutative,
one would not expect to get a groupoid model of it.  
  Instead we use Sieben's (unpublished) notion of {Fell bundles
over inverse semigroups}, which we believe should be thought of as
\stress{twisted \'etale groupoids with noncommutative unit space}.  

Recall that if $\G$ is an \'etale groupoid, then the set of open
bisections forms an inverse semigroup $\S$.  If for every $U$ in $\S$
we let $A_U$ be the set of elements in $C^*(\G)$ supported in $U$,
then $\A=\{A_U\}_{U\in\sS}$ is a Fell bundle over $\S$ whose
cross-sectional C*-algebra is isomorphic to $C^*(\G)$: although
Theorem (8.9) in \cite{\actions} is stated in a slightly different
language, this is precisely what it means.
Even if this result does not involve twists, in all likelihood it
may be proved for twisted groupoids as well.

Our second main result, Theorem \fcite{MainResult}{14.5}, shows that,
given a separable {\gc} pair $(A, B)$, there exists a
Fell bundle over a countable inverse semigroup $\S$ whose reduced
cross-sectional C*-algebra is isomorphic to $A$ in such a way that $B$
corresponds to the restriction of the bundle to the idempotent
semilattice of $\S$.

In a sense our result is not as complete as Renault's and its
predecessors since we go only as far as identifying the appropriate
Fell bundle.  In the commutative case (i.e.~when $B$ is commutative)
these results could perhaps be interpreted as further digging into
the structure of the Fell bundle and describing it by means of a
twisted groupoid.  Of course this task is not feasible when $B$ is
non-commutative, but perhaps there is more to say about our Fell
bundle than presently meets our eyes.

Let us now briefly discuss the methods used to prove our main
results.  While Renault's strategy is to first find the groupoid model
and then use it to prove the uniqueness of conditional expectations,
we do things in the opposite order. 

  Starting with a C*-algebra inclusion ``$B\subseteq A$", we say that
a \stress{slice} (Definition \fcite{DefineNormAndSlice}{10.1}) is a
closed linear subspace of the normalizer of $B$ in $A$ which is
invariant under left and right multiplication by elements of $B$.

Under suitable hypotheses every slice $M$ satisfies $M^*M\subseteq B$
so it may be thought of as Hilbert $B$-module, with inner product
$\<m,n\>=m^*n$.  In particular we may use Kasparov's stabilization
Theorem to prove the existence of a \stress{frame}, namely a sequence
$\{u_i\}_{i\in\N}$ of elements in $M$ such that
  $$
  m = \sum_{i\in \N} u_i\<u_i, m\>
  \for m\in M.
  $$
  Experimenting with variants of the notion of Jones-Watatani
index for conditional expectations \cite{\Jones,\Watatani}, we thought
of analyzing the series
  $$
  \tau = \sum_{i=1}^\infty\Cx(u_i)\Cx(u_i^*),
  $$
  where $\Cx$ is a conditional expectation onto $B$.  Surprisingly it
converges in the strict topology of the ideal
  $
  \rg M:= \closure{MM^*} \subseteq B,
  $
  as proved in Lemma \fcite{IntroduceTau}{11.5}.
  Since the Jones-Watatani index is central \scite{\Watatani}{1.2.8}
it is perhaps not surprising that $\tau$ lies in the center of the
multiplier algebra of $\rg M$. We moreover prove the following
curious and highly meaningful identity:
  $$
  \tau  mn^* = \Cx(m)\Cx(n^*)
  \for m,n\in M.
  $$
  One of its consequences is that $\|\Cx(m)\|= \|\tau\half
m\|$, for every $m$ in $M$, and hence the correspondence
  $$
  \Cx(m) \mapsto \tau\half m
  $$
  extends to an isometric map which turns out to be a {\vc} defined in
the ideal $\closure{\Cx(M)}$.  Assuming (\MaxPrime) we therefore deduce that
$\tau\half m\in B$, for every $m$, which puts us just a few steps away
from the proof of the uniqueness of conditional expectations.

The machinery developed to prove this result may then be used to
achieve the decomposition of a {\gc} pair $(A,B)$ via a
Fell bundle.  One proves without much difficulty that the set $\SAB$
of all slices forms an inverse semigroup (Proposition
\fcite{SABIsISG}{13.3}), and then the Fell bundle comes naturally by
setting $A_M=M$, for every $M\in\SAB$.

This paper is divided in two parts, the first one, comprising sections
(2--8), is designed not only to explain Sieben's theory of Fell
bundles over inverse semigroups (sections (2--3)), but also to define
the reduced cross-sectional C*-algebra and to prove that the fibers of
the bundle are faithfully represented.  The work in doing so is
considerable because, unlike the theory of Fell bundles over \stress{groups},
there is no readily available conditional expectation, a
useful tool for constructing \stress{regular
representations} in other contexts.

In order to overcome this difficulty we begin  in section (4) to study the much
simpler case of Fell bundles over \stress{semilatices}, i.e.~inverse semigroups
consisting of idempotent elements.
  In sections (5--8) we then develop a machinery designed to extend
states from the cross-sectional algebra of the bundle restricted to the
idempotent semilattice up to the whole cross-sectional algebra.  These
extended states are then used to get nontrivial representations.

Part two, comprising sections (9--14), is where we prove our main
results.  We begin with Section (9) by introducing and studying our
notion of {\vc s}. Theorem \fcite{ConcreteModelVC}{9.5}, one of the
main results of this section, gives a concrete model for these
gadgets, stating that every {\vc} is in fact represented by an element
which commutes with the subalgebra $B$, although it lives outside of
$A$.  Not all of the results proved in this section are used towards
our main results but they are included in the hope that the notion of
{\vc} might find applications elsewhere.

In section (10) we study slices and in section (11) we prove the
strict convergence of the series describing $\tau$ above and develop
its consequences.  Section (12) starts with our proposed
generalization of the notion of Cartan subalgebras and is where our
Theorem \fcite{UniqueCondexp}{12.3}
on uniqueness of conditional expectations is to be found.

In section (13) we construct the inverse semigroup and the Fell bundle
from a {\gc} pair which will be used in the next and final section to
provide our generalization (Theorem \fcite{MainResult}{14.5}) of
Renault's version of Feldman-Moore's Theorem.

  \font\bigrm=cmbxsc10 scaled 1440
  \font\medrm=cmcsc10 scaled 1440

\def\parte#1#2#3#4{\vfill\eject\null\vskip 0.5cm
  \hrule \bigskip
  \centerline{\bigrm #1}
  \bigskip\bigskip
  \centerline{\medrm #2}
  \def\a{#3}
  \def\b{#4}
  \ifx\a\b \else \medskip \centerline{\medrm #3}\fi
  \bigskip  \hrule
  \vskip 1cm}

\parte{PART ONE}
{Fell bundles over inverse semigroups and}
{their cross-sectional C*-algebras}{}

\section{Fell bundles over inverse semigroups}
  In this section we introduce the notion of Fell bundles over inverse
semigroups, which is the main object of interest in this work.  This
is essentially the same as Sieben's homonymous notion introduced in a
talk given in the Groupoid Fest, held at
Arizona State University in November 1998,
  with the title ``Fell bundles over inverse semigroups and r-discrete
groupoids".
  We thank Sieben for giving us access to his unpublished work which
we describe here with some small modifications.

  Throughout this section we let $S$ be an inverse semigroup, whose 
idempotent semilattice will be  denoted by $E(S)$.  We refer the
reader to \cite{\Lawson} for a comprehensive account of the theory of
inverse semigroups.

\definition \label DefineFellBundle A \stress{Fell bundle over $S$} is
a quadruple
  $$
  \A =  \Big(\{A_s\}_{s\in S}, \
  \{\mult_{s, t}\}_{s, t\in S}, \
  \{\starmap s\}_{s\in S}, \
  \{\incl_{t, s}\}_{s, t\in S,\, s\leq t}
  \Big)
  $$
  where, for each $s,t\in S$,
  \iaitem
  \aitem $A_s$ is a complex Banach space,
  \aitem $\mult_{s,t} : A_s\otimes A_t \to A_{st}$ is a linear map,
  \aitem $\starmap s : A_s \to A_{s^*}$ is a conjugate-linear
isometric map, and
  \aitem $\incl_{t, s}: A_s \hookrightarrow A_t$ is a linear isometric
  map for every $s\leq t$. 
  \medskip\noindent It is moreover required that for every $r,s,t\in
S$, and every $a\in A_r$, $b\in A_s$, and $c\in A_t$,
  \izitem
  \zitem $\mult_{rs, t}\big(\mult_{r, s}(a\otimes b) \otimes c\big) =
\mult_{r,st}\big(a\otimes\mult_{s, t}(b\otimes c) \big),$
  \zitem $\star{rs}{\mult_{r, s}(a\*b)} = \mult_{s^*, r^*}(\star
sb\*\star ra)$,
  \zitem $\star {s^*}{\star sa} = a$,
  \zitem $\|\mult_{r, s}(a\*b)\| \leq \|a\| \|b\|$,
  \zitem \zitemmark CstarIdentity $\|\mult_{r^*, r}\big(\star ra\*a\big)\| =
\|a\|^2$,
  \zitem \zitemmark Positivity $\mult_{r^*,r}\big(\star ra\*a\big) \geq 0$, in $B_{r^*r}$,
  \zitem \zitemmark Functorial if $r\leq s\leq t$, then $\incl_{t, r}=
\incl_{t, s}\circ \incl_{s,r}$,
  \zitem \zitemmark MultIndepend if $r\leq r'$, and $s\leq s'$, then
the diagrams
  \def\hup{\mult_{r, s}}%
  \def\hdn{\mult_{r', s'}}%
  \def\vl{\incl_{r', r}\otimes\incl_{s', s}}%
  \def\vr{\incl_{r's', rs}}%
  $$
  \hskip -50pt
  \matrix{
  \hfill A_r\otimes A_s \ & 
    \buildrel \hup \over \longrightarrow  &
    A_{rs} \hfill \vrule depth 15pt width 0pt \cr
  {\scriptstyle \vl} \downarrow \hskip 18pt && 
    \hskip 5pt \downarrow {\scriptstyle \vr}\cr
  \hfill  A_{r'}\otimes A_{s'} & 
    \buildrel \hdn \over \longrightarrow  &
    A_{r's'} \hfill \vrule height 18pt width 0pt 
  } 
  \hskip 25pt {\rm and}   \hskip 25pt 
  \matrix{
  \hfill A_s \ & 
    \buildrel \starmap s \over \longrightarrow  &
    A_{s^*} \hfill \vrule depth 15pt width 0pt \cr
  {\scriptstyle \incl_{s', s}} \downarrow \hskip 8pt && 
    \hskip 5pt \downarrow {\scriptstyle \incl_{s', s}}\cr
  \hfill  A_{s'} & 
    \buildrel \starmap{s'} \over \longrightarrow  &
    A_{s'^*} \hfill \vrule height 18pt width 0pt 
  }$$
  commute, 

\bigskip If no confusion is likely to arise we shall use the
simplified notations
  $$
  ab := \mult_{r, s}(a\*b)
  \and
  a^*:= \star r a,
  \eqmark SimplerNotation
  $$
  whenever $a\in A_r$ and $b\in A_s$.  Axioms
\lcite{\DefineFellBundle.i--vi} then take on the more familiar aspect:
  \def \zitemprime {\Item {{\rm(\rzitemplus')}}}
  \izitem
  \zitemprime $(ab)c = a(bc)$,
  \zitemprime $(ab)^* = b^*a^*$, 
  \zitemprime $a^{**} = a$,
  \zitemprime $\|ab\| \leq \|a\| \|b\|$,
  \zitemprime $\|aa^*\| = \|a\|^2$,
  \zitemprime $aa^* \geq 0$, in $B_{rr^*}$.

\bigskip If $s\leq t$, we shall often use the map $\incl_{t, s}$ to
identify $A_s$ as a subspace of $A_t$.  Axiom
\lcite{\DefineFellBundle.\MultIndepend} then says that the
multiplication and star operations are compatible with such an
identification.

Let us now discuss a few immediate consequences of the definition.

\state Proposition 
  \iaitem 
  \aitem If $e\in E(S)$, then $A_e$ is a C*-algebra.  Incidentally
this is the C*-algebra structure with respect to which
\lcite{\DefineFellBundle.\Positivity} refers to positivity.
  \aitem For every $s\in S$, one has that $\incl_{s, s}$ is the
identity map on $A_s$.
  \aitem If $e,f\in E(S)$, and $e\leq f$, then $\incl_{f, e}(A_e)$ is
a {\cltsi} in $A_f$.

\proof We skip (a), since it is obvious. Given $s\in S$, notice
that $\incl_{s, s}$ is an isometric linear map from $A_s$ to itself
which is idempotent by \lcite{\DefineFellBundle.\Functorial}.
Therefore $\incl_{s, s}$ must be the identity map on $A_s$, proving
(b).

With respect to (c), let $a\in A_e$ and $b\in A_f$.  Then
  $$
  \incl_{f,e}(a) b = 
  \mult_{f, f}\big(\incl_{f,e}(a) \*\incl_{f, f}(b)\big) 
  \={\lcite{\DefineFellBundle.\MultIndepend}} 
  \incl_{f, e}\big(\mult_{e, f}(a\*b)\big) \in \incl_{f, e}(A_e),
  $$
  and similarly $b\, \incl_{f,e}(a) \in \incl_{f, e}(A_e)$. This shows
that $\incl_{f, e}(A_e)$ is a two-sided ideal in $A_f$.  It is closed because
$A_e$ is a Banach space and $\incl_{f, e}$ is an isometric map.
\proofend

If $r, s, t\in S$ are such that $r, s\leq t$, let 
  $$
  v = ts^*sr^*r\  (\  = sr^*r = rs^*s).
  $$
  Then clearly $v\leq r,s$, and it is easy to see that $v$ is the
maximum among the elements of $S$ which are both smaller than $r$ and
$s$.  Therefore $v$ is effectively the \stress{infimum} of $r$ and $s$ with
respect to the natural order structure of $S$, so we shall denote $v$
by $r\wedge s.$ This is despite the fact that $S$ is \stress{not}
necessarily a semilattice with respect to its order: the existence of
$r\wedge s$ is only guaranteed once the set $\{r,s\}$ is known to admit an upper
bound, namely $t$.

\state Proposition
  Let $r, s, t\in S$ be such that $r, s\leq t$, and set $v = r\wedge s
= ts^*sr^*r$.  Then

  \begingroup \noindent \hfill \beginpicture
  \setcoordinatesystem units <0.0020truecm, -0.0015truecm> point at 0 0
  \put {$A_v$} at 0000 1000
  \put {$A_r$} at 1000 0000
  \put {$A_s$} at 1000 2000
  \put {$A_t$} at 2000 1000
  \arrow <0.15cm> [0.25,0.75] from 0200 0800 to 0800 0200
  \put {$\incl_{r, v}$} at 0300 300
  \arrow <0.15cm> [0.25,0.75] from 0200 1200 to 0800 1800
  \put {$\incl_{t, r}$} at 1700 300
  \arrow <0.15cm> [0.25,0.75] from 1200 0200 to 1800 0800
  \put {$\incl_{s, v}$} at 0350 1600
  \arrow <0.15cm> [0.25,0.75] from 1200 1800 to 1800 1200
  \put {$\incl_{t, s}$} at 1750 1600
  \endpicture \hfill\null \endgroup

  \bigskip\noindent is a pull-back diagram.

\proof Initially notice that the diagram commutes since 
  $$
  \incl_{t, r} \incl_{r, v} =   \incl_{t, v} = \incl_{t, s} \incl_{s, v}
  $$
  by \lcite{\DefineFellBundle.\Functorial}.  To prove the statement we
must show that whenever $a_r\in A_r$, and $a_s\in A_s$ are such that
  $$
  \incl_{t,r}(a_r)=  \incl_{t,s}(a_s),
  $$
  there exists a unique $a_v\in A_v$ such that 
  $$
  \incl_{r,v}(a_v) = a_r \and
  \incl_{s,v}(a_v) = a_s.
  $$
  Since all maps involved are injective, this is equivalent to saying
that 
  $$
  \incl_{t,v}(A_v) = \incl_{t,r}(A_r) \cap \incl_{t,s}(A_s).
  $$
  From
the commutativity of our diagram it easily follows that 
  $\incl_{t,v}(A_v) \subseteq \incl_{t,r}(A_r) \cap \incl_{t,s}(A_s)$,
so it suffices to prove the reverse inclusion, which we will do
shortly.

As already seen $A_{s^*s}$ is a C*-algebra, and the reader will have
no difficulty in proving that $A_s$ is a right Hilbert module over 
$A_{s^*s}$ relative to the operations
  $$
  x a := \mult_{s,s^*s}(x\*a) \and
  \<x,y\> := \mult_{s^*, s}\big(\star sx\*y\big),
  $$
  defined for all $a\in A_{s^*s}$, and $x,y\in A_s$.   Given an
approximate unit $\net u_i,i$ for $A_{s^*s}$, we then have by 
  \scite{\Jensen}{1.1.4} that
  $
  x=\lim_i x u_i,
  $
  for every $x\in A_s$.  

Given $y\in \incl_{t,r}(A_r) \cap \incl_{t,s}(A_s)$, write
  $$
  y =   \incl_{t,r}(a_r)=  \incl_{t,s}(a_s),
  $$
  where $a_r\in A_r$ and $a_s\in A_s$.  For every member $u_i$ of the
approximate unit above we have
  $$
  \incl_{t,s}(a_su_i) = 
  \incl_{t,s}\big(\mult_{s, s^*s}(a_s\*u_i)\big) \={(\DefineFellBundle.\MultIndepend)} 
  \mult_{t, t^*t}\big(\incl_{t, s}(a_s)\*\incl_{t^*t,
s^*s}(u_i)\big) \$=
  \mult_{t, t^*t}\big(\incl_{t, r}(a_r)\*\incl_{t^*t,
s^*s}(u_i)\big) \={(\DefineFellBundle.\MultIndepend)} 
  \incl_{t,rs^*s}\big(\mult_{r, s^*s}(a_r\*u_i)\big)  \in 
  \incl_{t,v}(A_v).
  $$
  So,
  $$
  y = \incl_{t,s}(a_s) = 
  \lim_i \incl_{t,s}(a_su_i) \in 
  \incl_{t,v}(A_v).
  \proofend
  $$

\section{Cross Sectional C*-algebra of a Fell Bundle}
  \def\clan{C^*\big(\LA /\NA\big)}%
  \def\iotaa{\iota_\A}%
  Throughout this section we fix a Fell bundle $\A$ over the inverse
semigroup $S$.

\definition \label NewDefineRep Let $B$ be a complex *-algebra.  A
\stress{pre-representation} of $\A$ in $B$ is a family $\Pi = \net
\pi_s,{s\in S}$, where for each $s\in S$
  $$
  \pi_s:A_s \to B
  $$
  is a linear map such that for all $s, t\in S$, all $a\in A_s$, and
all $b\in A_t$, one has
  \izitem
  \zitem 
  $
  \pi_{st}\big(\mu_{s, t}(a\* b)\big) =
  \pi_{s}(a)\pi_{t}(b),
  $
  \zitem $\pi_{s^*}(\star sa) = \pi_s(a)^*$.
  \medskip\noindent If moreover $\Pi$ satisfies 
  \zitem $\pi_t \circ \incl_{t, s} = \pi_s$, whenever $s\leq t$,
  \medskip\noindent we will say that $\Pi$ is a
\stress{representation}.

\bigskip
With the simplified notations given in \lcite{\SimplerNotation},
axioms \lcite{\NewDefineRep.i-ii} take on the simpler form:
  \izitem
  \zitemprime $\pi_{st}(ab) = \pi_s(a)\pi_{t}(b)$,
  \zitemprime $\pi_{s^*}(a^*) = \pi_s(a)^*$.

\medskip If $e\in E(S)$ we have already seen that $A_e$ is a
C*-algebra.  Given a pre-representation $\Pi$ of $\A$ in a C*-algebra
$B$, it is then immediate that $\pi_e$ is a *-homomorphism from $A_e$
to $B$.  In particular $\pi_e$ is necessarily
  contractive\fn{A linear map $T$ from a normed space $X$ to a normed
space $Y$ is \stress{contractive} if $\|T(x)\|\leq \|x\|$, for every
$x\in X$.}.
  Given any $s\in
S$, and $a\in A_s$, we then have that
  $$
  \|\pi_s(a)\|^2 = 
  \|\pi_s(a)^* \pi_s(a)\| = 
  \|\pi_{s^*}(\star sa)\pi_s(a) \|\$= 
  \|\pi_{s^*s}\big(\mu_{s^*, s}(\star sa\*a)\big)\| \leq 
  \|\mu_{s^*, s}(\star sa\*a)\big) \| \={(\DefineFellBundle.\CstarIdentity)}
  \|a\|^2.
  $$

With this we have shown:

\state Proposition \label RepContractive Given a pre-representation
$\Pi$ of the Fell bundle $\A$ in a C*-algebra $B$, one has that
  $$
  \|\pi_s(a)\| \leq  \|a\|
  \for s\in S \for a\in A_s.
  $$

\state Proposition \label IntersectProduct
  Let $\A$ be a Fell bundle over the inverse semigroup $S$ and let
$\Pi$ be a representation of $\A$ in the C*algebra $B$.  Then for
every idempotents $e,f\in E(S)$ one has that
  $$
  \pi_e(A_e) \cdot \pi_f(A_f) =
  \pi_{ef}(A_{ef}) =
  \pi_e(A_e)\cap \pi_f(A_f).
  $$

\proof Let $a\in A_e$ and $b\in A_f$.  Then, by
\lcite{\NewDefineRep.i} we have that 
  $$
  \pi_e(a) \pi_f(b) = \pi_{ef}\big(\mult_{e,f}(a\*b)\big) \in
\pi_{ef}(A_{ef}),
  $$
  proving that $\pi_e(A_e) \cdot \pi_f(A_f)
\subseteq\pi_{ef}(A_{ef}).$
  For any $c\in A_{ef}$ we have by \lcite{\NewDefineRep.iii} that 
  $$
  \pi_{ef}(c) = \pi_e\big(\incl_{e,ef}(c)\big) \in 
  \pi_e(A_e), 
  $$
  so 
  $\pi_{ef}(A_{ef})\subseteq\pi_e(A_e)$, and similarly
  $\pi_{ef}(A_{ef})\subseteq\pi_f(A_f)$, so 
  $$
  \pi_{ef}(A_{ef})\subseteq  \pi_e(A_e)\cap \pi_f(A_f).
  $$

To conclude we will prove that 
  $$
  \pi_e(A_e)\cap \pi_f(A_f)\subseteq \pi_e(A_e) \cdot \pi_f(A_f).
  $$
  It is evident that $\pi_e$ is a *-homomorphism from $A_e$ into $B$.
Therefore $\pi_e(A_e)$ is a closed *-subalgebra of $B$, the same
applying to $\pi_f(A_f)$.  If follows that $\pi_e(A_e)\cap\pi_f(A_f)$
is likewise a closed *-subalgebra of $B$.  

Given any $b\in \pi_e(A_e)\cap\pi_f(A_f)$ we may use 
  Cohen-Hewitt's factorization Theorem \scite{\CH}{32.22}
  to write $b = b_1b_2$, where $b_1,b_2\in \pi_e(A_e)\cap\pi_f(A_f)$,
so we see that
  $
  b\in \pi_e(A_e)\cdot\pi_f(A_f).
  $
  \proofend

\definition The \stress{cross-sectional C*-algebra of $\A$}, denoted
$C^*(A)$, is the universal C*-algebra generated by the disjoint union
  $$
  \dcup\limits_{s\in S} A_s,
  $$
  subject to the relations stating that the natural maps
  $$
  \piu_s:  A_s \to C^*(A)
  $$
  form a representation of $\A$ in $C^*(A)$.

From \lcite{\RepContractive} and \scite{\Blackadar}{1.2} one deduces the
existence of $C^*(\A)$ as well as its uniqueness up to isomorphism.
  For further reference let us spell out the universal property of
$C^*(\A)$ in detail:

\state Proposition \label UnivProp 
$C^*(\A)$ is a C*-algebra and $\Pi^u=\{\piu_s\}_{s\in S}$ is a
representation of $\A$ in $C^*(\A)$.  Moreover, 
given any representation $\Pi =
\net \pi_s,{s\in S}$ of the Fell bundle $\A$ in a C*-algebra $B$,
there exists a unique *-homomorphism $\Phi:C^*(\A)\to B$ such that
$\Phi\circ\piu_s = \pi_s$, for all $s\in S$.

One of our main goals is to show the existence of certain nontrivial
representations of $\A$, and also to show that each $\piu_s$ is
injective.
  In order to achieve this goal it is useful to have a more concrete
description of $C^*(\A)$, as follows.  Let
  $$
  \LA  = \bigoplus_{s\in S}A_s.
  $$
  For each $s\in S$ and each $a_s\in A_s$, denote by $a_s\d_s$ the
element of $\LA $ whose coordinates are all equal to zero except for the
$s^{\rm th}$ coordinate which is equal to $a_s$.  It is then clear
that for a generic element $a = (a_s)_s$ in $\LA $ one has
  $$
  a = \sum_{s\in S} a_s\d_s.
  $$
  In addition, this representation is clearly unique.

Define a multiplication and a star operation on $\LA $ in such a way
that for all $s,t\in S$, all $a\in A_s$, and all $b\in A_t$,
  $$
  (a\d_s)  (b\d_t) = \mult_{s,t}(a\*b)\d_{st}
  \and
  (a\d_s)^* = \star s{a}\d_{s^*}.
  $$  
  With the aid of \lcite{\DefineFellBundle.i--iii} one easily proves
that $\LA $ is a complex associative *-algebra.  Introducing a norm
on $\LA $ by
  $$
  \Big\|\sum_{s\in S} a_s\d_s\Big\| = 
  \sum_{s\in S} \|a_s\|,
  $$
  one readily checks that $\LA $ is a normed *-algebra.

\definition \label DefineRhoZero For each $s\in S$, let
  $$
  \pizero _s:a_s\in A_s \mapsto a_s\d_s\in \LA .
  $$

  It is easy to see that $\Pizero = \net \pizero _s, {s\in S}$ is a
pre-representation of $\A$ in $\LA $.  This is in fact a universal
pre-representation in the following sense:

\state Proposition \label StarRepVsPreRep Let $B$ be a *-algebra.
  If $\Pi = \net \pi_s,{s\in S}$ is a pre-representation of $\A$ in
$B$, then the map $\Phi:\LA \to B$ given by
  $$
  \Phi\Big(\sum_{s\in S}a_s\d_s\Big) = 
  \sum_{s\in S}\pi_s(a_s),
  $$
  is a *-homomorphism.  
  Conversely, 
  given any *-homomorphism $\Phi:\LA \to B$, consider for each $s\in S$,
the map   $\pi_s: A_s \to B$ given by
  $$
  \pi_s = \Phi\circ \pizero_s.
  $$
  Then $\Pi = \net \pi_s,{s\in S}$ is a pre-representation of $\A$ in
$B$.  
  In addition the correspondences
$\Pi\mapsto\Phi$, and $\Phi\mapsto\Pi$ described above are each other's
inverse, giving bijections between the set of all *-homomorphisms from
$\LA $ to $B$, and the set of all pre-representations of $\A$ in $B$.

\proof Left to the reader. \proofend

\state Proposition Let $B$ be a C*-algebra.  Then every *-homomorphism
$\Phi:\LA \to B$ is contractive.

\proof Let $\Pi$ be the pre-representation corresponding to $\Phi$
according to \lcite{\StarRepVsPreRep}. 
  Given any $a = \sum_{s\in S} a_s\d_s \in \LA $, we have that
  $$
  \|\Phi(a)\| =
  \Big\|\Phi\Big(\sum_{s\in S} a_s\d_s\Big)\Big\| =
  \Big\|\sum_{s\in S} \pi_s(a_s)\Big\| \leq
  \sum_{s\in S} \|\pi_s(a_s)\| \buildrel {(\RepContractive)} \over \leq
  \sum_{s\in S} \|a_s\| =
  \|a\|.
  \proofend
  $$

Observe that the pre-representation $\Pizero$ defined in
\lcite{\DefineRhoZero} is not necessarily a representation
since there is no reason why \lcite{\NewDefineRep.iii} holds.  In
order to force its validity we need to mod out certain elements of
$\LA $.

\state Proposition \label NewDefineIdealN Let $\NA$ be the linear
subspace of $\LA $ spanned by the set
  $$
  \Big\{a_s\d_s - \incl_{t, s}(a_s)\d_t :\
  s, t\in S,\ s\leq t,\  a_s\in A_s\Big\}.
  $$
  Then $\NA $ is a two-sided selfadjoint ideal of $\LA $.

\proof  Given $r,s,t\in S$, such that $s\leq t$, let $a_s\in A_s$, and
$b_r\in A_r$.  Then
  $$
  \Big(a_s\d_s - \incl_{t, s}(a_s)\d_t\Big)b_r\d_r =
  \mult_{s, r}(a_s\*b_r)\d_{sr} - \mult_{t, r}\big(\incl_{t,
s}(a_s)\*b_r\big)\d_{tr} 
  \={(\DefineFellBundle.\MultIndepend)} $$$$ =
  \mult_{s, r}(a_s\*b_r)\d_{sr} - \incl_{tr, sr}\big(\mult_{s,
r}\big(a_s\*b_r)\big)\d_{tr} \in \NA ,
  $$
  from where one sees that $\NA $ is a right ideal.  A similar reasoning
shows that $\NA $ is a left ideal as well.  Also
  $$
  \Big(a_s\d_s - \incl_{t, s}(a_s)\d_t\Big)^* =
  \star s{a}\d_{s^*} - \star t{\incl_{t, s}(a_s)}\d_{t^*} 
  \={(\DefineFellBundle.\MultIndepend)} $$$$ =
  \star s{a}\d_{s^*} - \incl_{t^*, s^*}\big(\star s{a_s}\big)\d_{t^*}
\in \NA ,
  $$
  so $\NA $ is selfadjoint.
  \proofend

The following result is self evident:

\state Proposition \label RepsVanishingOnNVsPhi
With respect to the correspondence $\Phi\leftrightarrow \Pi$ of
\lcite{\StarRepVsPreRep} one has that $\Phi$ vanishes on $\NA $ if and
only if $\Pi$ is a representation.

Observe that, by \lcite{\RepContractive} and \scite{\Blackadar}{1.2},
the enveloping C*-algebra of $\LA/\NA $, here denoted $\clan$,
exists.  It is our next main goal to prove that it is *-isomorphic to
$C^*(\A)$.

\definition \label DefineIotaA We denote by $\iotaa$ the composition
  \bigskip
  \begingroup \noindent \hfill \beginpicture
  \setcoordinatesystem units <0.0010truecm, 0.001truecm> point at 0 0
  \put {$\LA  \labelarrow q \LA/\NA \longrightarrow \clan,$} at 0000 0000
  \setquadratic
  \plot -2700 0400 -350 1000  2000 0400 / 
  \arrow <0.15cm> [0.25,0.75] from 2000 0400 to 2040 0380
  \put{$\iotaa$} at -300 1280
  \endpicture \hfill\null \endgroup

  \bigskip\noindent where $q$ is the quotient map and the unmarked
arrow is the canonical map from $\LA/\NA$ to its enveloping
C*-algebra.

\state Proposition
  For each $s\in S$, let $\piplus_s = \iotaa\circ \pizero _s$, where
$\pizero _s$ was defined in \lcite{\DefineRhoZero}, (see also Diagram
\fcite{DiagramaLegal}{3.14}) and let $\Piplus = \net \piplus_s, {s\in
S}$.  Then $\Piplus$ is a representation of $\A$ in $\clan$.

\proof Since $\Pizero $ is a pre-representation of $\A$ in $\LA $, it is
obvious that $\Piplus$ is a pre-representa\-tion of $\A$ in $\clan$.  To
prove that $\Piplus$ is a representation, let $s\leq t$, and let $a_s\in
A_s$.  Then by definition
  $a_s\d_s - \incl_{t, s}(a_s)\d_t\in \NA $, and hence 
  $$
  0 =
  \iotaa\big(a_s\d_s - \incl_{t, s}(a_s)\d_t\big) =
  \iotaa\big(\pizero _s(a_s) - \pizero _t(\incl_{t, s}(a_s)\d_t)\big) =
  \piplus_s(a_s) - \piplus_t(\incl_{t, s}(a_s)),
  $$
  proving that $\piplus_s = \piplus_t\circ\incl_{t, s}$, and hence
that $\Piplus$ is in fact a representation. \proofend

\state Proposition \label ClanIsomorphism There exists an isomorphism
$\Theta:\clan\to C^*(\A)$, such that $\Theta\circ\piplus_s = \piu_s$, for
every $s\in S$.

\proof
In order to prove the statement it is clearly enough to prove that
$\clan$ possesses the universal property described in
\lcite{\UnivProp} with respect to the representation $\Piplus$.  So we let
$\Pi = \net \pi_s,{s\in S}$ be any representation of $\A$ in a C*-algebra $B$.  
  Let $\Psi:\LA  \to B$ be given as in \lcite{\StarRepVsPreRep} in
terms of $\Pi$.  By \lcite{\RepsVanishingOnNVsPhi} we have that $\Psi$ 
vanishes on $\NA $ and hence it factors through $\LA/\NA $ giving a
*-homomorphism
  $\tilde\Psi: \LA/\NA  \to B$, 
  such that
  $$
  \tilde\Psi\big(q(a_s\d_s)\big) = \pi_s(a_s), 
  \subeqmark TildePsi
  $$
  whenever $a_s\in A_s$.  By the universality of the enveloping
C*-algebra, one has that $\tilde\Psi$ further factors through $\clan$, providing a
map $\Phi:\clan \to B$ such that the diagram

  \begingroup \noindent \hfill \beginpicture
  \setcoordinatesystem units <0.0010truecm, -0.001truecm> point at 0 0
  \put {$\LA/\NA $} at 1000 1000
  \arrow <0.15cm> [0.25,0.75] from 2200 1000 to 3400 1000
  \put {$\tilde\Psi$} at 2800 600
  \put {$B$} at 4000 1000
  \put {$\clan$} at 2300 3000
  \arrow <0.15cm> [0.25,0.75] from 1000 1500 to 1500 2500
  \arrow <0.15cm> [0.25,0.75] from 3100 2500 to 3700 1500
  \put {$\Phi$} at 3700 2200
  \endpicture \hfill\null \endgroup

\bigskip \noindent  commutes. In particular, for every $a\in A_s$ we
have 
  $$
  \pi_s(a_s) \={(\TildePsi)}
  \tilde\Psi\big(q(a_s\d_s)\big) =
  \Phi\big(\iotaa(a_s\d_s)\big) =
  \Phi\big(\iotaa\big(\pizero _s(a_s)\big)\big) =
  \Phi(\piplus_s(a_s)), 
  $$
  as desired.  That $\Phi$ is unique follows from the fact that
$\clan$ is generated by the union of the ranges of the $\piplus_s$.
\proofend

In view of the result above we shall henceforth identify $\clan$ and
$C^*(\A)$ bearing in mind that this identification caries $\piplus_s$
to $\piu_s$,  for every $s\in S$.
  The following diagram shows all relevant mappings:

  \bigskip
  \begingroup \noindent \hfill \beginpicture
  \setcoordinatesystem units <0.0010truecm, 0.0010truecm> point at 0 0
  \put {$\LA $} at 0000 0000
  \arrow <0.15cm> [0.25,0.75] from 0700 0000 to 2100 0000
  \put{$q$} at 1400 300
  \put{$\LA/\NA $} at 3000 000
  \arrow <0.15cm> [0.25,0.75] from 4000 0000 to 5100 0000
  \put{$\clan$} at 6500 000
  \arrow <0.15cm> [0.25,0.75] from 7900 0000 to 9200 0000
  \put{$\Theta$} at 8500 300
  \put{$C^*(\A)$} at 10000 000
  \setquadratic
  \plot 0200 0400 3100 1000  6000 0400 / 
  \arrow <0.15cm> [0.25,0.75] from 6000 0400 to 6040 0380
  \put{$\iotaa$} at 3100 1280
  \put{$A_s$} at 5000 -3000
  \arrow <0.15cm> [0.25,0.75] from 4500 -2800 to 0400 -400
  \put{$\pizero_s$} at 1800 -1900
  \arrow <0.15cm> [0.25,0.75] from 5200 -2600 to 6500 -400
  \put{$\piplus_s$} at 5300 -1300
  \arrow <0.15cm> [0.25,0.75] from 5500 -2800 to 9300 -400
  \put{$\piu_s$} at 8000 -1800
  \endpicture \hfill\null \endgroup

\bigskip
\centerline{\eightrm Diagram 
  \advseqnumbering\relax
  \current.
  \syslabel{DiagramaLegal}{\current}
  }

\section {Fell bundles over semilattices}
  \def\Ran{{\cal R}}%
  \def\Ranzero{\Ran_0}%
  Recall that a semilattice is a partially ordered set $E$ such that
for every $e,f\in E$ one has that there exists a larger element among
the members of $E$ which are smaller than both $e$ and $f$.  This
element is denoted by $e\wedge f$.  Viewed as a semigroup
under the operation ``$\wedge$",  one has that $E$ is an inverse
semigroup whose elements are idempotents.  Conversely any inverse
semigroup consisting of idempotents is obtained as above.

We will now study the special case of Fell bundles over semi-lattices.
So,  throughout this section we will let $\A$ be a Fell bundle over a
semi-lattice $E$.  Our main goal will be to construct a concrete
representation of $\A$ leading to a faithful representation of
$C^*(\A)$.

  \state Proposition \label BasicPropsRPi
  Given a representation $\Pi$ of $\A$ in a C*-algebra $B$, let 
  $$
  \Ranzero(\Pi) = \sum_{e\in E}\pi_e(A_e), 
  $$
  and let $\Ran(\Pi)$ be the closure of $\Ranzero(\Pi)$ within $B$.  Then
  \izitem
  \zitem $\Ranzero(\Pi)$ is a *-subalgebra of $B$,
  \zitem $\Ran(\Pi)$ is a closed *-subalgebra of $B$,
  \zitem $\pi_e(A_e)$ is a {\cltsi} of \/$\Ran(\Pi)$, for
every $e\in E$,
  \zitem $\Ranzero(\Pi)$ is a two-sided selfadjoint ideal of $\Ran(\Pi)$.

\proof  The first point follows easily from
\lcite{\NewDefineRep.i-ii}, and (ii) is an immediate consequence of (i).

In order to prove (iii) we must show that if  $a\in A_e$,
and $x\in \Ran(\Pi)$, then both $\pi_e(a)x$ and $x\pi_e(a)$ are in
$\pi_e(A_e)$.  Since $\pi_e(A_e)$ is closed, we may suppose without
loss of generality that $x=\pi_f(b)$, for some $f\in E$, and $b\in A_f$.
We have
  $$
  \pi_e(a)x = \pi_e(a) \pi_f(b) =
  \pi_{ef}\big(\mu_{e, f}(a\*b)\big) = \ldots
  $$
  Since $ef\leq f$, we have by \lcite{\NewDefineRep.iii} that the above
equals
 $$
  \ldots =
  \pi_e\Big( i_{e, ef}\big(\mu_{e, f}(a\*b)\big)\Big) \in
  \pi_e(A_e).
  $$
  That $x\pi_e(a)\in \pi_e(A_e)$ follows in the same way.  Finally,
(iv) holds as as consequence of the fact that the sum of ideals is an
ideal.
  \proofend

Before proceeding we need an elementary fact for which we were unable
to find a reference.

  \state Lemma \label AutoContinuous Let $A$ be a C*-algebra and let
$J$ be a (not necessarily closed) two-sided selfadjoint ideal of $A$.
Then any *-homomorphism $\Phi:J\to B$, where $B$ is a C*-algebra, is
contractive.

  \proof By replacing $A$ with its unitization we may assume that $A$
is unital.  We first claim that if $a\in J$, and $0\leq a\leq 1$, then
$\|\Phi(a)\|\leq 1$. For this let $b = a\sqrt{1-a^2\vrule height 9pt
width 0pt}$, so that $b\in J$.  Observe that
  $$
  a^4 + b^*b = a^2,
  $$
  hence 
  $$
  \Phi(a)^4 \leq
  \Phi(a)^4 + \Phi(b)^*\Phi(b) =
  \Phi(a)^2.
  $$
  Since $\Phi(a)$ is selfadjoint, one has that
  $$
  y:=\Phi(a)^2
  $$
  is a non-negative element of $B$ satisfying $y^2\leq y$, which
necessarily implies that $\|y\|\leq 1$.  Therefore
  $$
  \|\Phi(a)\|^2 = \|\Phi(a)^2\| = \|y\| \leq 1,
  $$
  proving our claim.

  If $a$ is any element of $J$, put $a':=\|a\|^{-2}a^*a$, so that $a'$
satisfies the hypotheses above and hence
  $$
  \|a\|^{-2}\|\Phi(a)\|^2 = 
  \|a\|^{-2}\|\Phi(a^*a)\| =  \|\Phi(a')\| \leq 1,
  $$
  from where the result follows.
  \proofend

We now present a very simple characterization of the
cross-sectional C*-algebra of $\A$.

\state Proposition \label IsomorphicIfOneToOne Let $\Pi$ be a
representation of $\A$ in a C*-algebra $B$ such that $\pi_e$ is
one-to-one for every $e$ in $E$.  Then $\Ran(\Pi)$ is *-isomorphic to
$C^*(\A)$.  In fact there exists a *-isomorphism
  $$
  \Phi:\Ran(\Pi)\to C^*(\A),
  $$
  such that $\Phi\circ\pi_e = \piu_e$,  for every $e\in E$.

\proof Let $\Pi'$ be any representation of $\A$ in a C*-algebra $B'$.
Given $x\in \Ranzero(\Pi')$, write
  $$
  x = \sum_e\pi'_e(a_e),
  $$
  where $a_e\in A_e$, for each $e\in E$, and the set $\{e\in E: a_e\neq 0\}$
is finite.  We begin by claiming that, if $f\in E$, and $b\in
A_f$, then
  $$
  x\pi'_f(b) = 
  \pi'_f \Big(  \sum_ei_{f, ef}\big(\mu_{e, f}(a_e\*b)\big)\Big).
  \subeqmark FormulaExpression
  $$
  To prove it we compute
  $$
  x\pi'_f(b) = 
  \sum_e\pi'_e(a_e)\pi'_f(b) =
  \sum_e\pi'_{ef}\big(\mu_{e, f}(a_e\*b)\big) \$= 
  \sum_e\pi'_f \Big(i_{f, ef}\big(\mu_{e, f}(a_e\*b)\big)\Big) =
  \pi'_f \Big(  \sum_ei_{f, ef}\big(\mu_{e, f}(a_e\*b)\big)\Big).
  $$ 
  Observe that, since $\Ran(\Pi')$ is generated by elements of the form
$\pi'_f(b)$, we have that
  $$
  x=0
  \quad \Longleftrightarrow \quad
  \big(\forall f\in E, \ \forall b\in A_f,\ x\pi'_f(b)=0\big).
  $$
  From \lcite{\FormulaExpression} we thus obtain the following
sufficient condition for $x$ to be zero:
  $$
  \sum_e\pi'_e(a_e)  = 0
  \quad  \Longleftarrow \quad
  \big(\forall f\in E, \ \forall b\in A_f,\ 
  \sum_ei_{f, ef}\big(\mu_{e, f}(a_e\*b)\big)=0
  \big).
  \subeqmark SufficCond
  $$
  Clearly all that has been said about $\Pi'$ holds, mutatis mutandis,
for $\Pi$.  Moreover, since we are supposing that $\pi_f$ is
one-to-one for every $f$, the sufficient condition \lcite{\SufficCond}
is also necessary in the case of $\Pi$.
  This implies that the expression 
  $$
  \Psi\Big(\sum_e\pi_e(a_e)\Big) = \sum_e\pi'_e(a_e)
  $$ 
  gives a well defined linear map from $\Ranzero(\Pi)$ onto 
$\Ranzero(\Pi')$, which can easily be proven to be a *-homomor\-phism as
well.  
  Since $\Ranzero(\Pi)$ is a two-sided selfadjoint ideal of $\Ran(\Pi)$ by
\lcite{\BasicPropsRPi.iv}, we deduce from \lcite{\AutoContinuous} that
$\Psi$ extends to give a *-homomorphism from $\Ran(\Pi)$ to $\Ran(\Pi')$,
which clearly satisfies 
  $$
  \Psi\circ\pi_e = \pi'_e
  \for e\in E.
  $$
  This shows that the conditions of \lcite{\UnivProp} are satisfied if
we replace $C^*(\A)$ and $\Pi^u$ by $\Ran(\Pi)$ and $\Pi$,
respectively.  Since the universal object is clearly unique up to
isomorphism the conclusion follows.
 \proofend

We thus see that, in order to obtain a concrete model for $C^*(\A)$, all
one needs is a faithful representation of $\A$.  In the remainder of
this section we shall obtain such a representation in a somewhat
canonical  way.

For each $e\in E$, consider the multiplier algebra $M(A_e)$ and let
$B$ be the subalgebra of
  $$
  \prod_{e\in E}M(A_e) 
  $$
  formed by the elements $m = (m_e)_{e\in E}$, satisfying 
  $$
  \|m \| := \sup_{e\in E} \|m_e\| < \infty.
  $$
  It is well known that $B$ is a C*-algebra with the norm defined
above.  Given $f\in E$ and $b\in A_f$ let, for every $e\in E$, 
  $$
  m^b_e = (L^b_e, R^b_e)\in M(A_e)
  $$
  be given my 
  $$
  L^b_e(a) =i_{e, fe}\big(\mu_{f,e}(b\*a)\big)
  \and
  R^b_e(a) =i_{e, ef}\big(\mu_{e,f}(a\*b)\big)
  \for a\in A_e.
  $$
  We then  let 
  $$
  \lambda_f(b) = (m^b_e)_{e\in E}\in B,
  $$
  leaving to the reader the routine verification that $\Lambda = \net
\lambda_e,{e\in E}$ is if fact a representation of $\A$.
  Observe that for every $b\in A_f$ we have that
  $$
  \|\lambda_f(b)\| \geq \|m^b_f\| = \|b\|,
  \eqmark LambdaIsometric
  $$
  because $m^b_f$ is the canonical multiplier of $A_f$ given by left
and right multiplication by $b$.   It follows that $\lambda_f$ is
one-to-one for every $f$.

As a consequence we deduce the main result of this section:

\state Corollary \label RegRepFaithfulForLatt Let $\A$ be a Fell
bundle over the semi-lattice $E$.  Then the cross-sectional C*-algebra
$C^*(\A)$ is *-isomorphic to $\Ran(\Lambda)$, where $\Lambda$ is the
representation of $\A$ constructed above.  In fact there exists a
*-isomorphism
  $$
  \Phi:\Ran(\Lambda)\to C^*(\A),
  $$
  such that $\Phi\circ\lambda_e = \piu_e$,  for every $e\in E$.

\proof Follows from \lcite{\IsomorphicIfOneToOne}. 
\proofend

In possession of a concrete representation we may now show that $\A$
is faithfully represented within $C^*(\A)$, at least in the case of
semi-lattices.

\state Corollary \label ResumoLattice Let $\A$ be a Fell bundle over
the semi-lattice $E$.  Then for every $e\in E$ one has that:
  \izitem
  \zitem $\piu_e$ is isometric,
  \zitem $\piu_e(A_e)$ is a {\cltsi} in $C^*(\A)$.

\proof The first point follows immediately from
\lcite{\LambdaIsometric} and \lcite{\RegRepFaithfulForLatt}.  The
second point is a consequence of \lcite{\BasicPropsRPi.iii} and
\lcite{\RegRepFaithfulForLatt}.
\proofend

\section{Support of linear functionals}
  In the above section we succeeded in obtaining, in a rather
elementary way, a nontrivial representation of $C^*(\A)$, in case $\A$
is a Fell bundle over a semi-lattice.  The case of Fell bundles over
general inverse semigroups is much more involving, requiring a careful
study of linear functionals defined on Fell bundles.  In the present
section we shall develop a few results about linear functionals on
C*-algebras which will be crucial in obtaining nontrivial
representations of general Fell bundles.

  The results presented here will most likely be well known by specialists
in the area and are included for completeness.

\state Proposition \label SupportConditions
  Let $A$ be a C*-algebra, let $J$ be a {\cltsi} of $A$, and let
$\phi$ be a {\clf} on $A$.  Then the following are equivalent:
  \izitem
  \zitem For every approximate unit $\net u_i,i$ for $J$, one has that
$\phi(a) = \lim_i\phi(au_i)$, for all $a\in A$.
  \zitem There exists a {\clf} $\psi$ on $J$ and an element $b$ in $J$
such that
  $
  \phi(a) = \psi(ab),
  $
  for all $a\in A$.
  \zitem For every approximate unit $\net u_i,i$ for $J$, one has that
$\phi(a) = \lim_i\phi(u_ia)$, for all $a\in A$.
  \zitem There exists a {\clf} $\psi$ on $J$ and
an element $b$ in $J$ such that
  $
  \phi(a) = \psi(ba),
  $
  for all $a\in A$.

\proof Applying 
  Cohen-Hewitt's factorization Theorem \scite{\CH}{32.22}
  to the restriction $\phi|_J$ we deduce that there exists a {\clf}
$\psi$ on $J$ and an element $b$ in $J$ such that
  $
  \phi(x) = \psi(xb),
  $
  for every $x\in J$.    

Assuming (i) choose an approximate unit $\net u_i,i$ for $J$ and
notice that for every $a\in A$ we have that
  $$
  \phi(a) = \lim_i\phi(au_i) = \lim_i\psi(au_ib) = \psi(ab),
  $$
  so (ii) follows.  
  Conversely, assuming that (ii) holds, let $\net u_i,i$ be any
approximate unit for $J$, and let $a\in A$.  One then has that
  $$
  \lim_i\phi(au_i) = 
  \lim_i\psi(au_ib) =
  \psi(ab) =
  \phi(a), 
  $$
  proving (ii).  In a similar way one proves that (iii) and
(iv) are equivalent.

We next prove that (ii) implies (iii).  For this let $\net u_i,i$ be
any approximate unit for $J$, and let $a\in A$.
  Then, assuming (ii) and observing that $ab\in J$, we have
  $$
  \lim_i\phi(u_ia)  = 
  \lim_i\psi(u_iab) =
  \psi(ab) = \phi(a),
  $$
  proving (iii).  Similarly (iv) implies (i), and the proof is
concluded.  \proofend

\definition If the equivalent conditions above hold we will say that
$\phi$ is \stress{supported} on $J$.

\state Proposition \label SupportedOnBigger Let $\phi$ be a {\clf} on the C*-algebra $A$, and let $I$ and $J$ be
{\cltsi s} of $A$ such that $I\subseteq J$.  If $\phi$ is supported on
$I$ then $\phi$ is also supported on $J$.

\proof
  By \lcite{\SupportConditions.ii} let $\psi$ be a {\clf} on $I$ and let $b\in I$, be such that
  $$
  \phi(a) = \psi(ab) \for a\in A.
  $$
  By Hahn-Banach's Theorem let $\chi$ be a continuous linear extension
of $\psi$ to $J$.  Then 
  $$
  \phi(a) = \chi(ab) \for a\in A,
  $$
  because $ab\in I$, so $\phi$ is supported on $J$.  
  \proofend

\state Proposition
  Let $A$ be a C*-algebra, let $J$ be a {\cltsi} of $A$,
and let $\phi$ be a state on $A$.  Denote by $\pi$ the GNS
representation of $A$ associated to $\phi$. Then the following are
equivalent:
  \izitem
  \zitem $\phi$ is supported on $J$.  
  \zitem The restriction of $\pi$ to $J$ is non-degenerated.

\proof Let $H$ be the space of $\pi$ and $\xi$ be the standard cyclic
vector.  Assuming (i) let $\net u_i,i$ be an approximate unit for $J$.
Then, for every $a, b\in A$ one has that
  $$
  \lim_i \<\pi(au_i)\xi,\pi(b)\xi\> =
  \lim_i \phi(b^*au_i) =
  \phi(b^*a) =
  \<\pi(a)\xi,\pi(b)\xi\>.
  $$
  Given that $\net \pi(au_i)\xi,i$ is a bounded net, we see that it
converges weakly to $\pi(a)\xi$.  This shows that
$\span\big(\pi(J)H\big)$ is weakly dense in $H$, and hence also
norm-dense.  Therefore $\pi|_J$ is non-degenerated, and (ii) follows.

Assuming that $\pi$ is non-degenerated, let once more $\net u_i,i$ be an
approximate unit for $J$.  It is well known that the net $\net
\pi(u_i),i$ converges strongly to the identity operator.  For every
$a\in A$ we therefore have that
  $$
  \phi(a) =
  \<\pi(a)\xi,\xi\> =
  \lim_i \<\pi(a)\pi(u_i)\xi,\xi\> =
  \lim_i \phi(au_i),
  $$  
  proving (i).  
  \proofend

\state Proposition \label Dicotomia Let $\phi$ be a pure state on a
C*-algebra $A$, and let $J$ be a {\cltsi} of $A$.  Then
either $\phi$ is supported on $J$ or $\phi$ vanishes on $J$.

\proof Let $\pi$ be the GNS representation of $A$ associated to $\pi$,
let $H$ be the space of $\pi$ and let $\xi$ be its associated cyclic
vector. If $\phi$ does not vanish on $J$ then obviously neither does
$\pi$.  Thus $\overline\span\big(\pi(J)H\big)$ is a nonzero closed
subspace of $H$, which is clearly invariant under $\pi$.  Given that
$\phi$ is pure we have that $\pi$ is irreducible, so
$\overline\span\big(\pi(J)H\big)=H$, and hence $\pi|_J$ is
non-degenerated.  So $\phi$ is supported on $J$. \proofend
  
\state Proposition \label SupportOnIntersect Let $\phi$ be a pure
state on a C*-algebra $A$, and let $I$ and $J$ be {\cltsi s} of $A$.
Suppose that $\phi$ is supported on both $I$ and $J$.  Then $\phi$ is
supported on $I\cap J$ as well. 

\proof Let $(\pi, \xi, H)$ be the GNS representation of $A$ associated
to $\phi$.  Choose approximate units $\net u_i,i$ and $\net v_j,j$ for
$I$ and $J$, respectively.
  As already mentioned both $\net \pi(u_i),i$ and $\net \pi(u_j),j$
converge strongly to the identity operator.  In particular
  $$
  \<\pi(u_{i_0})\xi, \xi\>\neq0,
  $$
  for some $i_0$.  Since
  $$
  \lim_j \phi(v_ju_{i_0}) =
  \lim_j \< \pi(v_j) \pi(u_{i_0})\xi,\xi\> =
  \<\pi(u_{i_0})\xi, \xi\> \neq 0 ,
  $$
  we have that
  $$
  \phi(v_{j_0}u_{i_0}) \neq0,
  $$
  for some $j_0$.  It follows that $\phi|_{I\cap J}$ is nonzero and
hence $\phi$ is supported on $I\cap J$ by \lcite{\Dicotomia}.
  \proofend

\section{Functionals on Fell bundles}
  Throughout this section we fix a Fell bundle $\A$ over the inverse
semigroup $S$.  Our main goal here will be to obtain a positive linear
functional on $\LA $, given a state on $A_e$, for some idempotent
$e\in E(S)$.  This will be the first step in obtaining representations
of $C^*(\A)$.

If $e\in E(S)$ and $s\in S$ are such that $e\leq s$, observe that the
map $\incl_{s, e}$ is an isometric embedding of $A_e$ into $A_s$.
Therefore we may view $A_e$ as a subspace of $A_s$.  Moreover, since
  $$
  e = se^*e =se \and e = ee^*s = es.
  $$
  we have that both $A_sA_e$ and $A_eA_s$ are contained in $A_e$.

\state Proposition \label extendingToAs 
  Let $e\in E(S)$ and $s\in S$ be such that $e\leq s$, and let $\phie$
be a {\clf} on $A_e$, then there exists a unique {\clf} $\tphie^s$ on
$A_s$, extending $\phie$, with $\|\tphie^s\| =\|\phie\|$, and such
that for every approximate unit $\net u_i,i$ for $A_e$ one has that
  $$
  \tphie^s(x) =
  \lim_i\phie(xu_i) =
  \lim_i\phie(u_ix) =
  \lim_i\phie(u_ixu_i)
  \for x\in A_s.
  $$

\proof By 
  Cohen-Hewitt's factorization Theorem \scite{\CH}{32.22}
  there exists a {\clf} $\psi$ on $A_e$ and an element $b$ in $A_e$
such that
  $$
  \phie(a) = \psi(ba)
  \for a\in A_e.
  $$
  Given $x\in A_s$ observe that $bx\in A_eA_s \subseteq A_e$, so it
makes sense to define
  $$
  \tphie^s(x) = \psi(bx)
  \for x\in A_s.
  $$
  Obviously $\tphie^s$ is then a {\clf} on $A_s$
extending $\phie$, whence $\|\tphie^s\| \geq \|\phie\|$.  Given an
approximate unit $\net u_i,i$ for $A_e$ one has, for every $x\in A_s$,
that
  $$
  \lim_i\phie(u_ix) =
  \lim_i\psi(bu_ix) =
  \psi(bx) =
  \tphie^s(x).
  $$
  On the other hand
  $$
  \lim_i\phie(xu_i) =
  \lim_i\psi(bxu_i) =
  \psi(bx) =
  \tphie^s(x),
  $$
  because $bx\in A_e$,  and hence $\net bxu_i,i$ converges to $bx$.
  To prove the last identity of the statement one uses 
  Cohen-Hewitt's Theorem once more to write
  $$
  \psi(a) = \chi(ac) \for a\in A_e,
  $$
  where $c\in A_e$, and $\chi$ is a {\clf} on
$A_e$. Then, for every $x\in A_s$ one has
  $$
  \lim_i\phie(u_ixu_i)  =
  \lim_i\psi(bu_ixu_i)  =
  \lim_i\chi(bu_ixu_ic) =
  \chi(bxc) =
  \psi(bx) = 
  \tphie^s(x).
  $$

  It remains to prove that $\|\tphie^s\| \leq \|\phie\|$.  For this
notice that for every $x\in A_s$
one has that
  $$
  |\tphie^s(x)| =
  \lim_i|\phie(xu_i)| \leq \limsup_i\|\phie\|\|x\|\|u_i\| \leq
  \|\phie\|\|x\|,
  $$
  concluding the proof. \proofend

\definition \label DefineExtendedFunct The functional $\tphie^s$
obtained above will be called the \stress{canonical extension} of
$\phie$ to $A_s$.

In order to facilitate computations with the canonical extension it is
useful to highlight the way it was constructed in the proof above:

\state Proposition \label WorkingCondition Under the hypothesis of
\lcite{\extendingToAs}, whenever
  $$
  \phie(a) = \psi(ba)
  \for a\in A_e,
  $$
  where $\psi$ is a {\clf} on
$A_e$, and $b\in A_e$, one has that 
  $$
  \tphie^s(x) = \psi(bx)
  \for x\in A_s.
  $$

\proof See the proof of \lcite{\extendingToAs}. \proofend

\definition Given $s, t\in S$, we will say that $s$ and $t$ are
  \izitem
  \zitem \stress{disjoint} if there is no $r\in S$ such that $r\leq
s,t$,
  \zitem \stress{$\A$-disjoint} if, for every $r\in S$ such that
$r\leq s,t$, one has that $A_r=\zerosp$.

Before proceeding we need the following elementary result which holds
for every inverse
semigroup:

\state Proposition \label SeEqualTe Let $s,t\in S$ and let $e\in
E(S)$.  Then the following are equivalent:
  \izitem
  \zitem $e\leq s^*t$, 
  \zitem $se=te$,  and $e\leq s^*s, t^*t$.

\proof
  (i) $\imply$ (ii): By (i) we have that $e=s^*te$.  So
  $$
  s^*se = s^*s(s^*te) = s^*te = e,
  $$
  proving that $e\leq s^*s$.  Since $e=e^*\leq (s^*t)^* = t^*s$, the
above reasoning gives $e\leq t^*t$.

  Proving that $se = te$ 
is equivalent to showing that $(se)^* = et^*$, which is to say that 
  \iaitem
  \aitem $se = se(et^*)se$,  and 
  \aitem $et^* = (et^*)se(et^*)$, 
  \medskip \noindent
  in view of the uniqueness of the adjoint in $S$.  Observing that
  $$
  et^*s =  (s^*te)^* = e^*=e,
  $$
  we have that the right-hand-side of (a) equals
  $$
  se(et^*)se =
  set^*se =
  see = se, 
  $$
  proving (a).
  To prove (b) notice that 
  $$
  (et^*)se(et^*) = 
  et^*set^* = 
  eet^* = et^*.
  $$

\noindent  (ii) $\imply$ (i):  This follows at once from 
  $s^*te = s^*se = e.$
  \proofend

\state Proposition \label PrimeiroVanishDisjoint Let $s$ and $t$ be
$\A$-disjoint elements in $S$ and let $e$ be an idempotent with $e\leq
s^*t$.  Given any {\clf} $\phie$ on $A_e$ one has
that
  $$
  \tphie^{s^*t}(a_s^*a_t)=0,
  $$
  for every $a_s\in A_s$, and $a_t\in A_t$.

\proof
By \lcite{\SeEqualTe} we have that $se = te$.  Denoting by $r=te$, it is then
clear that $r\leq s,t$, so $A_r=\{0\}$, by hypothesis.
  Next observe that by \lcite{\extendingToAs}
it is enough to show that
  $
  \phie(a_s^*a_tu)=0,
  $
  for every $u\in A_e$.  
  But since 
  $$
  a_tu\in A_tA_e \subseteq A_{te}=A_r=\{0\},
  $$
  we have that $a_tu=0$, and the statement follows.
  \proofend

\state Proposition \label DoubleExtend Let $e, s, t\in S$ be such that
$e$ is idempotent and $e\leq s \leq t$. Given a {\clf} $\phie$ on
$A_e$ one has that
  $$
  \tphie^s = \tphie^t\circ \incl_{t,s}.
  $$

\proof Pick a {\clf} $\psi$ on $A_e$,  and $b\in A_e$,  such that 
  $
  \phie(a) = \psi(ba)
  $
  for all $a\in A_e$.  Then, given any $x\in A_s$, we have
  $$
  \tphie^t\big(\incl_{t,s}(x)\big) \={(\WorkingCondition)} 
  \psi\big(b\incl_{t,s}(x)\big) = 
  \psi\big(\mult_{e, t}(b\*\incl_{t,s}(x))\big) 
  \={(\DefineFellBundle.\MultIndepend)} 
  \psi\big(\incl_{e, e}(\mult_{e, s}(b\*x))\big) \$=
  \psi(bx) = 
  \tphie^s(x). 
  \proofend
  $$

\state Proposition \label RestrictExtend Let $e,f\in E(S)$ and let
$s\in S$ be such that $e\leq f \leq s$. Given a {\clf} $\psi_f$ on
$A_{f}$ denote by $\phie$ its restriction to $A_{e}$. If $\psi_f$ is
supported on $A_{e}$ then the canonical extensions $\tphie^s$ and
$\tilde\psi_f^s$ coincide. 

\proof By \lcite{\SupportConditions.iv} let $\psi_1$ be a {\clf} on
$A_{e}$ and $b_1\in A_{e}$ be such that
  $$
  \psi_f(a) = \psi_1(b_1a)
  \for a\in A_{f}.
  $$
  By Hahn-Banach's Theorem we may suppose that $\psi_1$ is the
restriction to $A_{e}$ of a {\clf} $\psi_2$ on
$A_{f}$.
  Since 
  $
  \phie(a) = \psi_1(b_1a),
  $
  for all $a\in A_{e}$,
  we have by \lcite{\WorkingCondition} that 
  $$
  \tphie^s(x) = \psi_1(b_1x)
  \for x\in A_s.
  $$
  On the other hand, since 
  $
  \psi_f(a) = \psi_2(b_1a),
  $
  for all $a\in A_{f}$,
  we have, again by \lcite{\WorkingCondition}, that
  $$
  \tilde\psi_f^s(x) = \psi_2(b_1x)
  \for x\in A_s.
  $$
  However, since $b_1x\in A_{e}$, for all $x\in A_s$, we have that
  $$
  \tilde\psi_f^s(x) = \psi_2(b_1x) = \psi_1(b_1x) = \tphie^s(x).
  \proofend
  $$

\state Proposition \label PhiTildeIsPositive
  Let $e\in E(S)$, and let $\phie$ be a state on $A_e$.  For each
$s\in S$ such that $e\leq s$, denote by $\tphie^s$ the canonical
extension of $\phie$ to $A_s$, and let $\tphie$ be defined on $\LA $
by
  $$ 
  \tphie\Big(\sum_{s\in S}a_s\d_s\Big)=
  \sum_{s\geq e}\tphie^s(a_s).
  $$ 
  Then $\tphie$ is a positive {\clf} on $\LA $
with 
  $\|\tphie\|=\|\phie\|$.

\proof  It is clear that $\tphie$ is linear.  Moreover, if 
$\sum_{s\in S}a_s\d_s$ is a generic element of $\LA $, one has
that
  $$
  \Big|\tphie\Big(\sum_{s\in S}a_s\d_s\Big)\Big| \leq
  \sum_{s\geq e}|\tphie^s(a_s)| \leq 
  \sum_{s\geq e}\|\tphie^s\|\|a_s\| \={(\extendingToAs)}
  \|\phie\|  \sum_{s\geq e}\|a_s\| \leq
  \|\phie\|  \Big\|\sum_{s\in S}a_s\d_s\Big\|, 
  $$
  so $\|\tphie\|\leq \|\phie\|$.  Identifying $A_e$ with $A_e\d_e$ one
has that $\tphie$ extends $\phie$, so $\|\tphie\|=\|\phie\|$.
  In order to prove that $\tphie$ is positive, given $\ds a=\sum_{s\in
S}a_s\d_s\in \LA $, we need to show that $\tphie(a^*a)\geq0$.  We
have
  $$
  \tphie(a^*a) = 
  \tphie\Big(\sum_{s, t\in S}a_s^*a_t\d_{s^*t}\Big) =
  \sum_{s^*t\geq e}\tphie^{s^*t}(a_s^*a_t) \={(\SeEqualTe)}
  \sum_{
  \buildrel {\scriptstyle s^*s, t^*t \geq e} \over {\vrule height 6pt
width0pt se= te}
  }\tphie^{s^*t}(a_s^*a_t).
  \subeqmark SumAsStarAt
  $$
  Let 
  $
  X:= \{s\in S: s^*s\geq e\},
  $
  and consider the equivalence relation on $X$ defined by
  saying that 
  $s\sim t$, if and only if  $se=te$.
  Also let 
  $$
  X = \dcup_{\lambda\in \Lambda} X_\lambda,
  $$
  be the decomposition of $X$ into equivalence classes.
  Then, picking up from \lcite{\SumAsStarAt}, we have that
  $$
  \tphie(a^*a) = 
  \sum_{
  \buildrel {\scriptstyle s, t\in X} \over {s\sim t}
  }\tphie^{s^*t}(a_s^*a_t) =
  \sum_{\lambda \in \Lambda} \sum_{\ \ s,t\in X_\lambda} \tphie^{s^*t}(a_s^*a_t).
  $$
  To arrive at the conclusion it therefore suffices to show that
  $$
  \sum_{s,t\in X_\lambda} \tphie^{s^*t}(a_s^*a_t)\geq0
  \for \lambda\in\Lambda.
  \subeqmark SumOnLambda
  $$

Fixing $\lambda$ choose an approximate unit $\net u_i, i$ for $A_e$
and observe that $a_su_i\in A_{se} = A_{r_\lambda}$, where $r_\lambda$
is the common value of $se$ for all $s\in X_\lambda$.  We may then define
an element in $A_{r_\lambda}$ by the expression
  $$
  y_i=\sum_{s\in X_\lambda} a_su_i.
  $$
  Clearly 
  $
  y_i^*y_i\in 
  A_{r_\lambda^*r_\lambda} = 
  A_{(se)^*se} = 
  A_{es^*se} = 
  A_e,
  $ 
  so   we have by \lcite{\DefineFellBundle.\Positivity} that
  $$
  0\leq
  \phie (y_i^*y_i) =
  \phie\Big(\sum_{s,t\in X_\lambda} u_ia_s^*a_tu_i\Big) =
  \sum_{s,t\in X_\lambda} \phie(u_ia_s^*a_tu_i).
  $$
  Taking the limit as $i\to\infty$, we deduce that 
  $$
  0\leq
  \lim_i\sum_{s,t\in X_\lambda} \phie(u_ia_s^*a_tu_i) \={(\extendingToAs)}
  \sum_{s,t\in X_\lambda} \tphie^{s^*t}(a_s^*a_t),
  $$
  concluding the proof of \lcite{\SumOnLambda}.
  \proofend

\state Corollary \label SegundoVanishDisjoint
  Let $e\in E(S)$, and let $\phie$ be a state on $A_e$. If $s$ and $t$
are disjoint elements of $S$ then $\tphie$ vanishes on
$(A_s\d_s)^*(A_t\d_t)$.

\proof Given $a_s\in A_s$ and $a_t\in A_t$, one has that
  $$ 
  \tphie\big((a_s\d_s)^*(a_t\d_t)\big) = 
  \tphie\big(a_s^*a_t\d_{s^*t}).
  $$ 
  If $e\not\leq s^*t$, then the above vanishes by definition of
$\tphie$.  On the other hand, if $e\leq s^*t$, the above equals
  $ 
  \tphie^{s^*t}(a_s^*a_t),
  $ 
  which vanishes by \lcite{\PrimeiroVanishDisjoint}.
  \proofend

\section{Inducing pure states}
  \label PureStatesSection
  In the last section we have constructed states on $\LA $ from
states on each $A_e$.  These states do not necessarily vanish on the
ideal $\NA $ of \lcite{\NewDefineIdealN} and hence do not factor through
the algebra we are mostly interested in, namely $\LA/\NA $.  Here we
will improve the above construction in order to obtain states on
$\LA/\NA $ and hence nontrivial representations of $C^*(\A)$.

As before we fix a Fell bundle $\A$ over the inverse semigroup $S$.
Denote by $\E$ the restriction of $\A$ to the idempotent semilattice
$E(S)$, and 
recall that $\{\piu_e\}_{e\in E(S)}$ denotes the universal
representation of $\E$ in $C^*(\E)$. Incidentally, here we shall not use
the universal representation of $\A$.

By \lcite{\ResumoLattice.ii} we have that that $\piu_e(A_e)$ is a
{\cltsi} of $C^*(\E)$, for every $e\in E(S)$.

\definition Given a state $\phi$ on $C^*(\E)$, the \stress{support} of
$\phi$ is the set
  $$
  \supp(\phi) = \{e\in E(S): \phi \hbox{ is supported on } \piu_e(A_e)\}.
  $$

\state Proposition \label SupportIsFilter Let $\phi$ be a state on
$C^*(\E)$ and let $e, f\in E(S)$.
  \izitem \zitem
  If $e\in \supp(\phi)$, and $f\geq e$, then $f\in \supp(\phi)$.
  \zitem  If $\phi$ is pure and $e, f\in\supp(\phi)$, then $ef\in \supp(\phi)$.

\proof
  If $e\leq f$, then 
  $$
  \piu_f(A_f) \supseteq \piu_f\big(\incl_{f, e}(A_e)\big) =
\piu_e(A_e),
  $$
  so (i) follows from \lcite{\SupportedOnBigger}.
  As for (ii), it follows from \lcite{\SupportOnIntersect} and 
the fact that
  $
  \piu_{ef}(A_{ef}) =
  \piu_e(A_e)\cap \piu_f(A_f),
  $
  according to \lcite{\IntersectProduct}. \proofend

Recall that,  given two elements $e$ and $f$ of the  semilattice $E(S)$,
one says that $e\leq f$,  whenever $ef=e$.  For reasons that should
become clear let us introduce the \stress{reverse order} on $E(S)$ by
declaring that
  $$
  f \preceq e \iff e\leq f.
  $$
  In particular one has that 
  $$
  e,f\preceq ef
  \for e,f\in E(S).
  $$

  Given a pure state $\phi$ on $C^*(\E)$ observe that $\supp(\phi)$ is
a directed set under ``$\preceq$" by \lcite{\SupportIsFilter.ii}.
Therefore we may (and shortly will) use $\supp(\phi)$ as the set of
indices for a net.

Given $e\in \supp(\phi)$, the composition
  $$
  \phi_e := \phi\circ \piu_e 
  \eqmark DefinePhie
  $$
  is a state on $A_e$.  We therefore obtain by
\lcite{\PhiTildeIsPositive} a positive linear functional $\tphi_e$ on
$\LA $, and hence $\{\tphi_e\}_{e\in\supp(\phi)}$ is a net of
functionals on $\LA $.  This brings us to a main point.

\state Proposition \label BigExtension Given a pure state $\phi$ on
$C^*(\E)$, the net $\{\tphi_e\}_{e\in\supp(\phi)}$ constructed above
converges pointwise to a positive linear functional $\tphi$ on
$\LA $,  such that
  \izitem
  \zitem \zitemmark CalculoDeTPhi for every $s\in S$, and $a_s\in
A_s$, one has that
  $$
  \tphi(a_s\d_s) = \left\{\matrix{
  \tphi_e^s(a_s), & 
    \hbox{if there exists } e\in\supp(\phi), \hbox{ such that } e\leq s,
  \cr \cr 
  0, & \hbox{otherwise}, \hfill
  }\right.
  $$
  \zitem \zitemmark TphiExdendsMesmo for every $e\in E(S)$, and every
$a_e\in A_e$, one has that
  $
  \tphi(a_e\d_e) = \phi\big(\piu_e(a_e)\big),
  $
  \zitem \zitemmark TphiLessPhi $\|\tphi\| \leq \|\phi\|$,
  \zitem \zitemmark TphiVanisOnN $\tphi$ vanishes on the ideal $\NA $
defined in \lcite{\NewDefineIdealN}.

\proof To prove convergence of our net it is enough to show the
existence of 
  $$
  \lim_e\tphi_e(a_s\d_s)
  $$
  for every $s\in S$, and every $a_s\in A_s$. 

Suppose first that there does not exist $e\in\supp(\phi)$ such that
$e\leq s$.  By definition we have that $\tphi_e(a_s\d_s)=0$, for every
$e\in\supp(\phi)$, and hence 
  $$
  \lim_e\tphi_e(a_s\d_s) =0.
  \subeqmark LimitZero
  $$
  Suppose now that $e\leq s$, for some $e\in\supp(\phi)$.  We then
claim that for every $f\in\supp(\phi)$ with 
  $f\leq e$, 
  that is,  
  $f\succeq e$, 
  one has that
  $$ 
  \tphi_f(a_s\d_s)=\tphi_{e}(a_s\d_s).
  $$ 
  To prove it observe that
  $$
  \phi_f = \phi\circ \piu_f =
  \phi\circ \piu_{e} \circ \incl_{e, f} =
  \phi_{e} \circ \incl_{e, f},
  $$
  so, identifying $A_f$ as an ideal of $A_{e}$, we have that
$\phi_f$ coincides with the restriction of $\phi_{e}$ to $A_f$.  It
is also evident that $\phi_{e}$ is supported on $A_f$, so we have by 
\lcite{\RestrictExtend} that $\tphi_f^s = \tphi_{e}^s$.  Therefore
  $$
  \tphi_f(a_s\d_s) = \tphi_f^s(a_s) =
  \tphi_{e}^s(a_s) =
  \tphi_{e}(a_s\d_s),
  $$
  proving our claim, and hence that 
  $$
  \lim_f\tphi_f(a_s\d_s) = \tphi_{e}(a_s\d_s).
  \subeqmark LimitNonZero
  $$
  This concludes the proof of the convergence of our net.
  Letting $\tphi$ denote the pointwise limit of the $\tphi_e$, we have
that \lcite{\CalculoDeTPhi} follows from \lcite{\LimitZero} and
\lcite{\LimitNonZero}.
  Since each $\tphi_e$ is positive by \lcite{\PhiTildeIsPositive} it
is clear that $\tphi$ is also positive.

In order to prove \lcite{\TphiExdendsMesmo}, suppose first that
$e\in\supp(\phi)$.  Then, using \lcite{\CalculoDeTPhi} and 
observing that the extension of $\phi_e$ to $A_e$ is $\phi_e$ itself,
we have that
  $$
  \tphi(a_e\d_e) =
  \tphi_e^e(a_e) = 
  \phi_e(a_e) = 
  \phi\big(\piu_e(a_e)\big).  
  $$
  Suppose next that $e\notin\supp(\phi)$.  This  means that  $\phi$ is
not supported on $\piu_e(A_e)$, which implies that $\phi$ vanishes on
$\piu_e(A_e)$ by \lcite{\Dicotomia}, and hence the right-hand-side of
\lcite{\TphiExdendsMesmo} equals zero.
  Because of \lcite{\SupportIsFilter.i}, no member of the support of
$\phi$ is dominated by $e$, so $\tphi(a_e\d_e)=0$, by
\lcite{\CalculoDeTPhi}.  This concludes the proof of \lcite{\TphiExdendsMesmo}.

Addressing \lcite{\TphiLessPhi}
notice that
for every $e\in
\supp(\phi)$, one has that
  $$
  \|\tphi_e\| \={(\PhiTildeIsPositive)}
  \|\phi_e\| = 
  \|\phi\piu_e\| \leq
  \|\phi\|\|\piu_e\| \buildrel{(\RepContractive)} \over \leq
  \|\phi\|.
  $$
  Taking the limit one gets $\|\tphi\| \leq \|\phi\|$.

  Let us now prove that $\tphi$ vanishes on $\NA $,  so let
$s,t\in S$ be such that $s\leq t$, and let $a_s\in A_s$.  By
definition of $\NA $ \lcite{\NewDefineIdealN}, our task is then to prove
that
  $$
  \tphi\big(a_s\d_s - \incl_{t, s}(a_s)\d_t\big) = 0.
  \subeqmark OurTask
  $$
  Consider the following statements:
  \iaitem
  \aitem There exists $e\in \supp(\phi)$ such that $e\leq s$.
  \aitem There exists $e\in \supp(\phi)$ such that $e\leq t$.

  \medskip Since $s\leq t$, it is clear that (a) implies (b).  It is
thus impossible for (a) to be true and (b) false.   All other cases
will be treated separately.
  \def\case#1{\medskip\noindent {\tensc Case #1}:}%

  \case 1 (a) is false and (b) is true. 

Choose $e_0\in\supp(\phi)$ such that $e_0\leq t$. Then, for every
$e\leq e_0$ (i.e. $e\succeq e_0$), we have that 
  $$
  se =
  ts^*se =
  tes^*s =
  es^*s,
  $$
  so $se$ is idempotent, and clearly $se\leq s$.  Since we are
assuming (a) to be false, we deduce that $se$ is not in
$\supp(\phi)$ and hence that $\phi$ is not supported in
$\piu_{se}(A_{se})$.  So $\phi$ vanishes on $\piu_{se}(A_{se})$ by
\lcite{\Dicotomia}, and   in particular
  $$
  \phi_e\incl_{e, se} = 
  \phi  \piu_e\incl_{e, se} = 
  \phi  \piu_{se} = 0.
  \subeqmark PhiESEZero
  $$
  Still assuming that $e\leq e_0$, we next claim that 
  $$
  \tphi_e^t\big(\incl_{t, s}(a_s)) = 0.
  $$
  Fixing an approximate unit $\net u_i,i$ for $A_e$ we have that
  $$
  \tphi_e^t\big(\incl_{t, s}(a_s)) \={(\extendingToAs)}
  \lim_i\phi_e\big(\incl_{t, s}(a_s)u_i\big) =
  \lim_i\phi_e\Big(\mult_{t, e}\big(\incl_{t, s}(a_s)\*u_i\big)\Big) 
  \={(\DefineFellBundle.\MultIndepend)} $$$$ =
  \lim_i\phi_e\Big(\incl_{e, se}\big(\mult_{s, e}(a_s\*u_i)\big)\Big)
  \={(\PhiESEZero)} 0, 
  $$
  proving our claim.
  It follows that
  $$
  \tphi_e\big(a_s\d_s - \incl_{t, s}(a_s)\d_t\big) \={(\PhiTildeIsPositive)}
  \tphi_e^t\big(\incl_{t, s}(a_s)) = 0,
  $$
  for every $e\succeq e_0$, and hence \lcite{\OurTask} is proved.

\case 2 (a) and (b) are both false.

For every $e\in \supp(\phi)$  we have that 
  $
  \tphi_e\big(a_s\d_s - \incl_{t, s}(a_s)\d_t\big) = 0, 
  $
  by \lcite{\PhiTildeIsPositive}, so \lcite{\OurTask} follows
immediately.

\case 3 (a) and (b) are both true.

  Choose $e_0\in \supp(\phi)$ such that $e_0\leq s$.  Then, for all
$e\leq e_0$ (i.e. $e\succeq e_0$), we have that 
  $$
  \tphi_e\big(a_s\d_s - \incl_{t, s}(a_s)\d_t\big) = 
  \tphi_e^s(a_s) - \tphi_e^t\big(\incl_{t, s}(a_s)\big)
  \={(\DoubleExtend)} 0,
  $$
  and hence \lcite{\OurTask} follows. \proofend

It is interesting to remark that the convergence above takes place
even if $\C$ is given the discrete topology!

\state Proposition \label TerceiroVanishDisjoint
  Let $\phi$ be a pure state on $C^*(\E)$.  If $s$ and $t$ are
disjoint elements of $S$ then $\tphi$ vanishes on
$(A_s\d_s)^*(A_t\d_t)$.

\proof Since $\tphi$ is the limit of the $\tphi_e$,  the result follows
from \lcite{\SegundoVanishDisjoint}. \proofend

\section{Representations of $C^*(\A)$ and the reduced cross-sectional algebra}
  As before we fix a Fell bundle $\A$ over the inverse semigroup $S$
and let $\E$ be the restriction of $\A$ to the idempotent semilattice
$E(S)$.

  Given a pure state $\phi$ on $C^*(\E)$ one has that the functional $\tphi$
provided by \lcite{\BigExtension} is continuous and hence extends
to the Banach *-algebra completion of $\LA$.
  Using \scite{\Fell}{VI.19.3 and VI.19.5} we deduce that $\tphi$
\stress{generates} a *-representation $\RiLA$ of (the completion of $\LA$ and
hence also of) $\LA$ on a Hilbert space $H_{\tphi}$.
  To be precise, $H_{\tphi}$ is the Hilbert space obtained as the
Hausdorff\fn{Meaning that one must first mod out vectors of norm
zero.}  completion of $\LA$ under the pre-inner-product
  $$
  \<x,y\>_{\tphi} = \tphi(y^*x)
  \for x,y\in \LA.
  \eqmark InnerProductViaphi
  $$
  The representation $\RiLA$ itself is given by
  $$
  \RiLA(x)\hat y = \widehat{xy}
  \for x,y\in\LA,
  \eqmark DefineRiLA
  $$
  where the ``hat" notation indicates the canonical map 
$y\in \LA \mapsto \hat y\in H_{\tphi}$, given by the completion process.

\definition 
We shall refer to $\RiLA$ as the \stress{GNS representation of $\LA $
associated
to $\tphi$}.

Notice however that, in the absence of a bounded approximate unit for
$\LA$, the existence of a cyclic vector for $\RiLA$ does not follow
from the standard arguments.

Regardless of cyclic vectors we have the following criterion for the
vanishing of $\RiLA(x)$, for every given $x\in \LA $:
  $$
  \RiLA(x)=0 \ \Longleftrightarrow \ \forall y,z\in \LA ,  \
\tphi(z^*xy)=0.
  \eqmark CriterionForVanishRho
  $$

\definition \label DefineReduced By the \stress{reduced C*-algebra of
$\A$} we shall mean the C*-algebra $\CstarRed(\A)$ obtained as the
Hausdorff completion of $\LA $ under the C*-seminorm
  $$
  \tnorm x = \sup_\phi\|\RiLA(x)\|,
  $$
  where $\phi$ runs in the set of all pure states of $C^*(\E)$.

\state Proposition There exists a surjective map $\Lambda:C^*(\A) \to
\CstarRed(\A)$, such that the diagram
  \bigskip
  \begingroup \noindent \hfill \beginpicture
  \setcoordinatesystem units <0.0015truecm, -0.0015truecm> point at 0 0
  \put {$\LA $} at 000 000
  \put {$C^*(\A)$}       at 2500 -1000
  \put {$\CstarRed(\A)$} at 2500  1000
  \put {$\Lambda$}   at 2800 -0100
  \put {$\iotaa$}    at 1000 -1000
  \put {$\iotaaRed$} at 1000  1000
  \arrow <0.15cm> [0.25,0.75] from 2500 -600 to 2500 550
  \arrow <0.15cm> [0.25,0.75] from 600 -350 to 1700 -850
  \arrow <0.15cm> [0.25,0.75] from 600  350 to 1700  850
  \endpicture \hfill\null \endgroup

\medskip\noindent commutes, where $\iotaaRed$ is the canonical map arising from
the completion process.

\proof If $\phi$ is a pure state of $C^*(\E)$ recall that $\tphi$
vanishes on $\NA $ by \lcite{\BigExtension.\TphiVanisOnN}.  So it follows
from \lcite{\CriterionForVanishRho} that $\RiLA$ also vanishes on $\NA $,
and hence $\RiLA$ factors through $\LA/\NA $, and hence also through
its enveloping C*-algebra,  namely $C^*(\A)$.
  It follows that 
  $
  \|\RiLA(x)\| \leq \|\iotaa(x)\|,
  $
  and hence also that
  $$
  \|\iotaaRed(x)\| = 
  \tnorm x \leq \|\iotaa(x)\|
  \for x\in \LA ,
  $$
  from where the conclusion follows. \proofend

\definition \label DefinePiur For each $s\in S$ we shall let
$\piur_s:A_s \to \CstarRed(\A)$ be the composition of maps
  $$
  A_s \labelarrow{\pizero_s}\LA  \labelarrow{\iotaaRed } \CstarRed(\A).
  $$

Notice that 
  $$
  \piur_s= 
  \iotaaRed \circ \pizero_s =
  \Lambda \circ \iotaa \circ \pizero_s =
  \Lambda\circ \piu_s,
  $$
  and hence $\Pi^\reduced = \{\piur_s\}_{s\in S}$ is a representation
of $\A$ in $\CstarRed(\A)$.

We may now finally prove a non-triviality result relating to $C^*(\A)$.

\state Lemma
  Let $n\geq 1$, let $s_1, \ldots, s_n\in S$, be pairwise
$\A$-disjoint elements, and let $a_i\in A_{s_i}$, for each $i=1, 2,
\ldots, n$.  If
  $$
  \sum_{i=1}^n \piur_{s_i}(a_i) = 0,
  $$
  then $a_i=0$, for every $i$.

\proof
  Supposing by contradiction that some $a_k\neq0$, let $e_k=
s_k^*s_k$.  By \lcite{\DefineFellBundle.\CstarIdentity} and
\lcite{\ResumoLattice.i} we have that $\piu_{e_k}(a_k^*a_k)\neq0$,
so we may choose a pure state $\phi$ on $C^*(\E)$ such that 
$\phi\big(\piu_{e_k}(a_k^*a_k)\big)\neq0$.

Obviously $\phi$ does not vanish on the ideal $\piu_{e_k}(A_{e_k})$
and consequently $\phi$ is supported there by \lcite{\Dicotomia}.  In
other words, $e_k\in\supp(\phi)$.  Since
  $$
  0 =
  \sum_{i=1}^n \piur_{s_i}(a_i) =
  \iotaaRed \Big(\sum_{i=1}^n a_i\d_{s_i}\Big),
  $$
  we must have that $\RiLA\Big(\sum_{i=1}^n a_i\d_{s_i}\Big)=0$.
  Given any $u\in A_{e_k}$ we obtain from
\lcite{\CriterionForVanishRho} that 
  \def\mixterm#1#2{(a_{#1}u\d_{s_{#1}e_k})^*(a_{#2}u\d_{s_{#2}e_k})}%
  $$
  0 =
  \tphi\Big( \big(u\d_{e_k}\big)^*\big(\sum_{i=1}^n a_i\d_{s_i}\big)^*
\big(\sum_{i=1}^n a_i\d_{s_i}\big)\big(u\d_{e_k}\big)
\Big) = 
  \tphi\Big( \sum_{i, j=1}^n \mixterm ij \Big).
  \subeqmark TPhiDeuZero
  $$
  By hypothesis we have that $s_i$ and $s_j$ are $\A$-disjoint for
$i\neq j$.  The same is therefore also the case for $s_ie_k$ and
$s_je_k$, so we have by \lcite{\TerceiroVanishDisjoint} that 
$\tphi\big(\mixterm ij\big)=0$.  This means that the cross
terms in \lcite{\TPhiDeuZero} all vanish and we are left with 
  $$
  0 =
  \tphi\Big( \sum_{i=1}^n \mixterm ii\Big) \geq
  \tphi\big(\mixterm kk\big) \$=
  \tphi(u^*a_k^*a_ku\d_{e_k})  \={(\BigExtension.i)}
  \tphi_{e_k}^{e_k}(u^*a_k^*a_ku)  =
  \phi_{e_k}(u^*a_k^*a_ku) \={(\DefinePhie)}
  \phi\big(\piu_{e_k}(u^*a_k^*a_ku)\big).
  $$
  Letting $u$ run through an approximate unit for $A_{e_k}$ we deduce
that $\phi\big(\piu_{e_k}(a_k^*a_k)\big)=0$, in contradiction to the
choice of $\phi$.  Therefore $a_k$ must be zero and the proof is
concluded.
  \proofend

The following is an immediate consequence:

\state Corollary \label MapsAreFiel 
  The maps 
  $$
  \piu_s: A_s \to C^*(\A),
  $$ and $$
  \piur_r: A_s \to \CstarRed(\A),
  $$
  are injective for every $s\in S$.

\state Corollary \label IdentifyEinA There are monomorphisms
  $
  \Phi: C^*(\E) \to C^*(\A)
  $, 
  and 
  $
  \Phi_\reduced:  C^*(\E) \to \CstarRed(\A),
  $
  such that for every $e\in E(S)$ the diagram 
  \medskip
  \begingroup \noindent \hfill \beginpicture
  \setcoordinatesystem units <-0.0010truecm, -0.0010truecm> point at 0 0
  \put {$\CstarRed(\A)$} at -2700 0000
  \put {$C^*(\E)$} at 0000 0000
  \put {$C^*(\A)$} at 2700 0000
  \put {$A_e$} at 0000 2000
  \arrow <0.15cm> [0.25,0.75] from -700 0000 to -1800 0000 
    \put{$\Phi_\reduced$} at -1200 -300
  \arrow <0.15cm> [0.25,0.75] from 800 0000 to 1800 0000
    \put{$\Phi$} at 1200 -300
  \arrow <0.15cm> [0.25,0.75] from 0000 1600 to 0000 0400
    \put{$\piu_e$} at 0300 900
  \arrow <0.15cm> [0.25,0.75] from  0400 1700 to  2300 0400
    \put{$\piu_e$} at 1800 1200
  \arrow <0.15cm> [0.25,0.75] from -0400 1700 to -2300 0400
    \put{$\piur_e$} at -1800 1200
  \endpicture \hfill\null \endgroup

\noindent  commutes. 

\proof Follows immediately from \lcite{\MapsAreFiel} and
\lcite{\IsomorphicIfOneToOne}. \proofend

Notice that in the statement of \lcite{\IdentifyEinA} we are denoting
by $\piu_e$ the canonical maps in two different contexts, namely that
of the Fell bundle $\A$ and that of the restricted bundle $\E$.
Nevertheless by this result we may view $C^*(\E)$ as a subalgebra of
$C^*(\A)$ and, once this identification is made, the two meanings of
$\piu_e$ are reconciled.

In a forthcoming paper we plan to further develop the theory of
cross-sectional algebras, but for the time being we present the
following simple fact:

\state Proposition \label ApproxIdentInCstarE Any (bounded)
approximate unit for $C^*(\E)$ is an approximate unit for both
$C^*(\A)$ and $\CstarRed(\A)$.

\proof Let $\{\app_i\}_i$ be an approximate unit for $C^*(\E)$.  Given
$s\in S$, it follows from
  \scite{\Jensen}{1.1.4} that $\piu_s(A_s)$ is the
closed linear span of $\piu_s(A_s)\piu_s(A_s)^*\piu_s(A_s)$.
Observing that
  $$
  \piu_s(A_s)\piu_s(A_s)^* \subseteq \piu_{ss^*}(A_{ss^*}) \subseteq C^*(\E),
  $$
  we see that $C^*(\A)$ equals the closed linear span of
$C^*(\E)C^*(\A)$, and therefore also of $C^*(\E)C^*(\A)C^*(\E)$.  

Incidentally, this is to say that $C^*(\A)$ is the {\cltsi} generated
by $C^*(\E)$.  The result for $C^*(\A)$ then follows easily and the
case of $\CstarRed(\A)$ may be treated similarly. \proofend

\parte{\bigrm PART TWO}{{\Gc} subalgebras}{}{}

\section{Virtual commutants}
  Having followed the first few steps into the general theory of Fell
bundles over inverse semigroups, we now wish to show that they appear
naturally in certain situations.  In fact we wish to show that, under
suitable hypothesis, the inclusion of a closed *-subalgebra $B$ of a
C*-algebra $A$ is, in all respects, the same as the inclusion
  $$
  C^*(\E) \subseteq \CstarRed (\A),
  $$
  for a Fell bundle $\A$.  For a while we will forget about Fell
bundles and will concentrate instead on inclusion of C*-algebras
satisfying, to begin with, the following:

\sysstate{Standing Hypothesis}{\rm}{\label StandingOne 
  From now on we assume that $A$ is a \stress{separable}
C*-algebra and $B\subseteq A$ is
a closed *-subalgebra containing an approximate unit 
  $\{\app_i\}_i$ for $A$.}

\definition \label VirtuComm A \stress{\vc} of $B$ in $A$ is an
$A$-valued linear map $\phi$ defined on a {\cltsi} $J$ of $B$, such
that
  \izitem 
  \zitem  $\phi(bx) = b\phi(x)$, and 
  \zitem  $\phi(xb) = \phi(x)b$,
  \medskip\noindent for all $x\in J$, and $b\in B$.  Should  we want
to highlight the domain of $\phi$, we will write $(J,\phi)$ in place
of $\phi$.

Notice that, under the conditions above, both $J$ and $A$ are
$B$-bimodules.  Thus conditions (i) and (ii) simply say that $\phi$ is
a $B$-bimodule map.

\state Proposition \label ElementaryVirtual Let $(J,\phi)$ be a
{\vc} of $B$ in $A$.  Then
the range of $\phi$ is contained in $JAJ$.

\proof By 
  Cohen-Hewitt's factorization Theorem \scite{\CH}{32.22} 
  any $x\in J$ may be written as $x=yzw$, with $y,z, w\in J$.  Then
  $$
  \phi(x) = \phi(yzw) = y\phi(z)w \in JAJ.
  \proofend
  $$

\state Proposition Every {\vc} $(J, \phi)$ is bounded.

\proof Employing the closed graph Theorem, we must prove that if
$\{x_n\}_{n\in\N} \subseteq J$ is such that $x_n\to 0$ and
$\phi(x_n)\to a\in A$, then $a=0$.  

Notice that, since $\phi(x_n)\in JA$, by
\lcite{\ElementaryVirtual}, we have that $a = \lim_n\phi(x_n)\in
\closure{JA}$.  One may then prove that $a = \lim_iu_ia$, for every
(bounded) approximate unit $\{u_i\}_i$ for $J$.

In order to prove that $a=0$, it is therefore enough to prove that
$ya=0$, for every $y\in J$.  We have
  $$
  ya =
  \lim_{n\to\infty} y\phi(x_n) =
  \lim_{n\to\infty} \phi(yx_n) =
  \lim_{n\to\infty} \phi(y)x_n = 
  \phi(y)  \lim_{n\to\infty} x_n = 0.
  \proofend
  $$

\bigskip
Here is a procedure for obtaining a fairly general example of a
{\vc}.  Suppose that $A$ is a closed *-subalgebra of some
other C*-algebra $C$, that is $B\subseteq A \subseteq C$, and that
$\T$ is an element of $C$ which commutes with $B$.
  Letting $J$ be the subset of $B$ given by
  $$
  J=\{b\in B: \T b \in A\},
  $$
  it is easy to see that $J$ is a {\cltsi} of
$B$. Defining 
  $$
  \phi(x)=\T x \for x\in J,
  $$
  one then has that $(J,\phi)$ is a {\vc} of $B$ in $A$.

The following result essentially states that the above example is the
most general one.

\state Theorem \label ConcreteModelVC Let $\phi$ be a {\vc} of $B$ in
$A$.  Supposing that $A$ is faithfully represented on a Hilbert space
$H$ (in which case we shall identify $A$ with its image within
$B(H)$), there exists $\T\in B(H)$ such that for every $x\in J$,
  \izitem
  \zitem $\T x\in A$,
  \zitem $\phi(x) = \T x$,  and 
  \zitem $\T$ commutes with $B$.

\proof 
  Observe that by 
  Cohen-Hewitt's factorization Theorem \scite{\CH}{32.22},
  every element $\xi\in\closure{JH}$ can be written as $\xi=x\eta$,
for $x\in J$, and $\eta\in H$.  In particular $\closure{JH} = JH$.

Letting $\{u_i\}_i$ be a bounded approximate unit for $J$, we claim that the
net $\{\phi(u_i)\xi\}_i$ converges to some element in $JH$, for every
$\xi\in JH$.  In fact, writing $\xi=x\eta$, as above, we have
  $$
  \phi(u_i)\xi = 
  \phi(u_i)x\eta = 
  \phi(u_ix)\eta \buildrel i\to\infty \over \longrightarrow
  \phi(x)\eta.
  $$
  That $\phi(x)\eta$ lies in $JH$ follows from
\lcite{\ElementaryVirtual}.  The correspondence 
  \ $\xi \mapsto\lim_i\phi(u_i)\xi$ \
  therefore defines a bounded linear map $\T:JH\to JH$ such that 
  $$
  \T(x\eta) = \phi(x)\eta
  \for x\in J \for \eta\in H.
  $$
  Incidentally this shows that $\T$ does not depend on the choice of
the approximate unit above.

Declaring $\T$ to be zero on $JH^\perp$ we get an
extension of $\T$ to $H$ which, by abuse of language, will still be 
denoted by $\T$.  

We next claim that $\T$ commutes with $B$.  Notice that since $J$ is
an ideal in $B$, the space $JH$ is invariant under $B$,  and clearly
also under $\T$.
So, in order to prove that $\T b=b\T$, for any given $b\in B$, it is
enough to prove that
  $
  \T b(\xi) = b\T(\xi),
  $
  for all $\xi\in JH$  or, equivalently, that 
  $$
  \T b(x\xi) = b\T (x\xi)
  \for x\in J \for \xi\in H.
  $$
  We have
  $$
  \T b(x\xi) =
  \phi(bx)\xi =
  b\phi(x)\xi  =
  b\T(x\xi).
  $$
  This proves (iii).
  We next prove (ii).
  Given $\zeta\in H$, write $\zeta = y\xi+\eta$, with $y\in
J$, $\xi\in H$, and $\eta \in JH^\perp$.  Using the fact that $\T$
vanishes on $JH^\perp$, we have
  $$
  \T x(\zeta) = 
  \T x(y\xi+\eta)  =
  \T (xy\xi) = 
  \phi(xy)\xi = 
  \phi(x)y\xi = 
  \phi(x)(y\xi+\eta)  =
  \phi(x)(\zeta),
  $$
  where the penultimate equality, namely that $\phi(x)\eta=0$, follows
from the facts that $\phi(x)\in AJ$, by \lcite{\ElementaryVirtual},
and that $J$ vanishes on $JH^\perp$.  So $\phi(x) = \T\xi$, taking
care of (ii), and hence also of (i).
  \proofend

\definition We shall say that $B$ satisfies property (\MaxPrime) 
if, for any
{\vc} $\phi$ of $B$ in $A$, one has that the range of $\phi$ is
contained in $B$.

As usual we will use the notation $B'\cap A$ to refer to the relative
commutant of $B$ in $A$, namely 
  $$
  B'\cap A = \{a\in A: ab=ba, \forall b\in B\}.
  \label DefineRelativeCommutant
  $$

\state Proposition \label CompareclassicalNotion If $B$ satisfies
property (\MaxPrime) then $B'\cap A\subseteq B$.

  \proof Let $a\in B'\cap A$.  Defining
  $$
  \phi:x\in B \mapsto xa \in A,
  $$
  it is easy to see that 
$(B,\phi)$ is a {\vc}.  So the hypothesis implies that the
range of $\phi$ lies in $B$ or, equivalently, that $Ba\subseteq B$.
  Under \lcite{\StandingOne} we deduce that $a\in B$, hence proving
that $B'\cap A\subseteq B$.  \proofend

The converse of the above theorem is not necessarily true.  To
describe a counter-example let $B$ be a (necessarily non-unital)
C*-algebra whose center reduces to $\{0\}$, such as the algebra of
compact operators on an infinite dimensional Hilbert space.  Let $X$
be any compact Hausdorff topological space with more than one point
such as the real interval $[0,1]$, and let us view $B$ as the
subalgebra $1\*B\subseteq C(X)\*B =C(X, B)$ formed by the constant
functions.  It is then easy to see that
  $
  B'\cap C(X,B) = \{0\} \subseteq B,
  $
  but if one chooses a non-constant function $f\in C(X)$, the map
  $$
  \phi: b\in B\mapsto f\*b \in C(X)\*B,
  $$
  gives a {\vc} whose range is not contained in $B$.

This phenomena however does not occur in the realm of abelian
algebras.

\state Proposition \label GeneralizeMaxAbel Assume that $B$ is
abelian.  Then $B$ satisfies property (\MaxPrime) if and only if
$B'\cap A\subseteq B$, which incidentally is to say that $B$ is
maximal abelian.

\proof    The implication ``$\Rightarrow$"  follows immediately from
\lcite{\CompareclassicalNotion}.  Conversely,
let $(J, \phi)$ be a
{\vc} of $B$ in $A$.  Then, for every $x\in J$ and $b\in
B$, one has that
  $$
  b\phi(x) = 
  \phi(bx) = 
  \phi(xb) = 
  \phi(x)b,
  $$
  so $\phi(x)\in B'\cap A\subseteq B$, proving that $B$ satisfies
property (\MaxPrime).
  \proofend
  
\section{Slices and the normalizer of $B$ in $A$}
  We will now begin to study the notion of normalizer in the context
of C*-algebras.  This notion was introduced by Kumjian \cite{\Kumjian}
to study abelian subalgebras and was also used by Renault
\cite{\Renault}.

  In this section, as well as throughout the rest of this work, we keep
\lcite{\StandingOne} in force.

\definition \label DefineNormAndSlice The \stress{normalizer} of $B$
in $A$ is the subset
  $$
  N(B) = \{a\in A: a^*Ba\subseteq B,\ aBa^*\subseteq B\}.
  $$
  A \stress{slice} is any closed linear subspace $M\subseteq N(B)$,
such that both $BM$ and $MB$ are contained in $M$.

Observe that $N(B)$ is a closed set.  In addition it is closed under
multiplication and under adjoint, but it is not necessarily closed
under addition.  On the other hand, by definition a slice is required
to be closed under addition, but it is not necessarily closed under
multiplication or adjoint.

\state Proposition \label MMStarContainedInB Let $M$ be a slice.  Then
the sets $M^*BM$, $MBM^*$, $M^*M$ and $MM^*$ are contained in $B$.

\proof Let $m,n\in M$.  Then,  since $M$ is a linear subspace  we have
that $m+i^kn\in M$, for $k=0,1,2,3$.  Therefore, given $b\in B$,
  $$
  m^*bn={1\over 4}\sum_{k=0}^3i^{-k}(m+i^kn)^*b(m+i^kn) \in B,
  $$
  hence $M^*BM\subseteq B$.  Applying the same reasoning to the slice
$M^*$ we deduce that $MBM^*\subseteq B$.  

Still assuming that $m,n\in M$, let $\{\app_i\}_i$ be as in
\lcite{\StandingOne}.  Then
  $$
  m^*n = \lim_i m^*\app_in \in B,
  $$
  so $M^*M\subseteq B$, and similarly $MM^*\subseteq B$.
  \proofend

\state Corollary \label SlicesTRO Let $M$ be a slice. Then
  \izitem 
  \zitem $MM^*M\subseteq M$, which is to say that $M$ is a
\stress{ternary ring of operators} (cf. {\rm \cite{\Zettl}}),
  \zitem $M^*M$ and $MM^*$ (linear span of products) are two-sided
self-adjoint ideals of $B$.

\proof For (i) it is enough to notice that
  $MM^*M \subseteq BM \subseteq M$.  Point (ii) is trivial.
  \proofend

The following notation will be useful:

\definition Given a slice $M$ we say that
  \izitem
  \zitem the \stress{source} of $M$, denoted $\so M$, is the {\cltsi}
$\overline{M^*M}$ (closed linear span) of $B$, and 
  \zitem the \stress{range} of $M$, denoted $\rg M$, is the {\cltsi}
$\overline{MM^*}$ of $B$.

\state Proposition \label SliceForEveryone Every element $a\in N(B)$
lies in some slice.

\proof We claim that $M :=\overline{BaB}$ (closed linear span) is a
slice. It is clear that $BM,MB\subseteq M$, and in order to prove that
$M$ is contained in $N(B)$, it is enough to notice that $MBM^*$ and
$M^*BM$ are contained in $B$.  Finally observe that by
\lcite{\StandingOne} one has that
  $a = \lim_i \app_ia\app_i \in M$.  \proofend

  \section{Frames and conditional expectations}
  We will now discuss the behavior of slices under conditional
expectations.
We must however start by discussing the notion of frames, with which
we begin the present section.

  We continue to work under \lcite{\StandingOne},  observing that it is only
from now on that we make use of the separability of $A$.

\state Proposition \label GotFrame Let $M$ be a slice.  Then there
exists a countable family $\{u_i\}_{i\in\N}$ of elements of $M$ such that
  $$
  m = \sum_{i\in \N} u_iu_i^*m
  \for m\in M,
  $$
  where the sum converges unconditionally.  Any such family will be
called a \stress{frame}.

\proof We regard $M$ as a right Hilbert module over $\Bu$ (unitization
of $B$) under the
inner product
  $$
  \<m,n\> = m^*n \for m,n\in M.
  $$
  Since $A$ is separable it follows that $M$ is also separable and
hence, by Kasparov's stabilization Theorem 
  \scite{\Kasparov}{3.2}, \scite{\Jensen}{1.1.24},
  there exists a unitary operator
  $$
  U: \ell_2(\Bu) \to M\oplus \ell_2(\Bu).
  $$
  Let $\{e_i\}_{i\in\N}$ be the canonical basis of $\ell_2(\Bu)$ and
let $u_i=P\big(U(e_i)\big)$, where $P$ denotes the orthogonal pro\-jec\-tion
from $M\oplus \ell_2(\Bu)$ onto $M$.
  Given any $\xi\in \ell_2(\Bu)$ it is easy to see that 
  $$
  \xi = \sum_{i\in\N}e_i\<e_i,\xi\>,
  $$
  where the sum converges unconditionally.
  In particular, given $m\in M$,  we have for $\xi=U\inv(m)$, that
  $$
  m = P(m) = P\big(U(\xi)\big) =
  \sum_{i\in\N}P\big(U(e_i)\big)\<e_i,\xi\> =
  \sum_{i\in\N}u_i\<e_i,\xi\>.
  $$
  We also have that
  $$
  \<e_i,\xi\> =
  \<U(e_i),U(\xi)\> =
  \<U(e_i),m\> =
  \<U(e_i),P(m)\> \$=
  \<P\big(U(e_i)\big),m\> =
  \<u_i,m\> =
  u_i^*m,
  $$
  which combines with the calculation above to give the result.
  \proofend

\state Proposition \label MyNiceApproximateUnit Let $M$ be a slice and
let $\{u_i\}_{i\in\N}$ be a frame for $M$.  Also let $\F$ be the
set of all finite subsets of\/ $\N$, ordered by inclusion and viewed
as a directed set.  Then
  \izitem 
  \zitem $\Big\{\sum_{i\in F}u_iu_i^*\Big\}_{F\in\F}$ is an
approximate unit for $\rg M$,
  \zitem  $\|\sum_{i\in F}u_iu_i^*\|\leq 1$,  for every $F\in\F$.

  \proof For every $F\in\F$, let 
  \def\psum{s}%
  $
  \psum_F=\sum_{i\in F}u_iu_i^*,
  $
  and if $k\in\N$, let $\psum_k = \psum_{\{1,2,\ldots, k\}}$. Consider
the bounded operator
  $$
  T_k:m\in M \mapsto \psum_km\in M.
  $$
  By hypothesis we have that $T_k$ converges pointwise to the identity
operator on $M$ and hence $\{T_k\}_k$ is a uniformly bounded set by
the Banach-Steinhaus Theorem.  By \scite{\tpa}{4.7} we have that 
$\|\psum_k\| =\|T_k\|$, whence
  $$
  \sup_k\|\psum_k\| =   \sup_k\|T_k\| <\infty.
  $$  
  Given $F\in \F$, pick $k\in\N$ such that $F\subseteq\{1,2, \ldots,k\}$,
and observe that
  $$
  0\leq \psum_F = \sum_{i\in F}u_iu_i^* \leq \sum_{i\leq k}u_iu_i^* = \psum_k,
  $$
  so $\ds\sup_{F\in\F}\|\psum_F\| <\infty$, as well.
If $m,n\in M$ we have that
  $$
  \lim_{F\to\infty} \psum_F mn^* =   \sum_{i\in\N} u_iu_i^*mn^* = mn^*,
  $$
  so $\lim \psum_F b=b$, for every $b\in MM^*$.  Since the $\psum_F$ are
bounded, the last identity holds for every $b$ in the closed linear
span of
$MM^*$, thus proving (i).  With respect to (ii), it follows from the
elementary result immediately below.
  \proofend

\state Proposition Let $\{v_i\}_i$ be an increasing approximate identity
of a C*-algebra $A$ consisting of positive elements. Then $\|v_i\|\leq 1$,
for all $i$.

\proof For  every $i\leq j$ and every selfadjoint $a$ in $A$ one has that
  $
  a\half v_ia\half  \leq a\half v_ja\half .
  $
  Taking the limit on $j$ we conclude that 
  $$
  a\half v_ia\half  \leq a.
  $$
  With $a=v_i$ we get 
  $
  v_i^2\leq v_i,
  $ 
  so $\|v_i\|\leq 1$.
\proofend

We now begin our study of conditional expectations.

\state Proposition
  \label FirstLocationOfEM
  Let $\Cx:A\to B$ be a conditional expectation, and let $M$ be a
slice.  Then $\Cx(M) \subseteq \rg M\cap \so M$.

\proof Recall from
  \scite{\Jensen}{1.1.4} that if $\{v_i\}_i$ is an
approximate identity for $\rg M$, then $\lim_iv_im = m$, for all
$m\in M$.  Therefore
  $$
  \Cx(m) =
  \lim_i \Cx(v_im) =
  \lim_i v_i\Cx(m) \in \rg M.
  $$
  Similarly $\Cx(m^*) \in \rg {M^*} = \so M$, whence 
  $$
  \Cx(m) = \Cx(m^*)^* \in \big(\so M\big)^* = \so M.
  \proofend
  $$

The following is a key technical result:

\state Lemma
  \label IntroduceTau
  Let $\Cx:A\to B$ be a conditional expectation, and let $M$ be a slice.
Then there exists a unique central positive element $\tau$ in the multiplier
algebra of $\rg M$ such that $\|\tau\|\leq1$, and 
  $$
  \tau  mn^* = \Cx(m)\Cx(n^*)
  \for m,n\in M.
  \subeqmark IntroduceTauEquation
  $$
  Moreover, if $\{u_i\}_{i\in\N}$ is any frame for $M$, one has that 
  $\tau = \sum_{i\in\N}\Cx(u_i)\Cx(u_i^*)$,
  the series being unconditionally convergent in the strict topology.

\proof Since $\rg M$ is the closure of $MM^*$, it is obvious that
$\tau$ is unique.  As for existence let $\{u_i\}_{i\in\N}$ be a frame
for $M$. We claim that the series
  $$
  \sum_{i\in\N}\Cx(u_i)\Cx(u_i^*)
  \subeqmark TheSeriesToConvStrictly
  $$
  converges unconditionally in the strict topology of the multiplier
algebra of $\rg M$.  This means that the net $\{\tau_F\}_{F\in\F}$,
where $\F$ was defined in \lcite{\MyNiceApproximateUnit}, and
  $$
  \tau_F = 
  \sum_{i\in F}\Cx(u_i)\Cx(u_i^*)
  \for F\in\F,
  $$
  converges strictly.  In turn this is to say that for every $b\in\rg
M$, the nets $\{b\tau_F\}_{F\in\F}$, and $\{\tau_Fb\}_{F\in\F}$
converge in the norm topology of $\rg M$.  We treat first the case in
which $b=mn^*$, with $m,n\in M$.  Given $F\in \F$,  we have
  $$
  \tau_Fmn^* =
  \sum_{i\in F}\Cx(u_i)\Cx(u_i^*)mn^* =
  \sum_{i\in F}\Cx(u_i)\Cx(u_i^*mn^*) =
  \sum_{i\in F}\Cx(u_i)u_i^*m\Cx(n^*) \$=
  \sum_{i\in F}\Cx(u_iu_i^*m)\Cx(n^*) =
  \Cx\Big(\sum_{i\in F}u_iu_i^*m\Big)\Cx(n^*).
  $$
  By \lcite{\GotFrame} we then conclude that 
  $$
  \lim_{F\to\infty} \tau_Fmn^* =
  \Cx(m)\Cx(n^*) =
  \lim_{F\to\infty} mn^* \tau_F,
  \subeqmark BothLimits
  $$
  where the second equality is actually a consequence of the first, by
taking adjoints.  To prove convergence in case $b$ is not necessarily
of the form $mn^*$, it is now enough to show that
$\{\tau_F\}_{F\in\F}$ is a bounded net, but this may be proved with
the help of 
  \scite{\TakeOne}{III.3.4.iii} 
  as follows:
  $$
  0\leq \tau_F = \sum_{i\in F}\Cx(u_i)\Cx(u_i^*) \leq
  \sum_{i\in F}\Cx(u_iu_i^*)  = 
  \sum_{i\in F}u_iu_i^* \buildrel (\MyNiceApproximateUnit.ii) \over \leq 1.
  $$
  This shows that \lcite{\TheSeriesToConvStrictly} does converge as
indicated, so we let $\tau$ be its sum.  Obviously $\|\tau\|\leq 1$.
Identity \lcite{\IntroduceTauEquation} then follows immediately from
\lcite{\BothLimits}, and so does the fact that $\tau$ is
central. \proofend  

We will soon be interested in the space $\closure{\tau\half M}$, when we will
need the following two elementary results:

\state Lemma \label ElementaryOnFunctions Let $a$ be a positive real
number.  Given $\alpha,\beta, \varepsilon>0$, there exists a
continuous scalar valued function $g$ on $[0,a]$ such that $g(0)=0$,
and
  $$
  \big|t^\beta - g(t)t^\alpha\big|< \varepsilon
  \for t\in [0,a].
  $$

  \proof
  If $\beta>\alpha$ it is enough to take $g(t) = t^{\beta-\alpha}$,
so we suppose that $\beta\leq\alpha$.  For every positive integer $n$, 
let
  \bigskip  \bigskip
  $$
  \matrix{
  \hbox to 15mm {\hfill}
  &
  g_n(t) = \left\{\matrix{ 
    n^{\alpha-\beta+1}t \  , & \hbox{ if } t\leq {1\over n}, \cr\cr
    \hfill t^{\beta-\alpha} \hfill \ , & \hbox{ if } t\geq {1\over n}.}
  \right.
  &
  \hbox to 20mm {\hfill}
  &
  
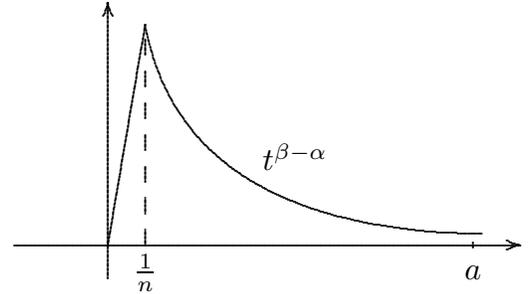
\beginpicture
  \setcoordinatesystem units <0.00125truecm, 0.00075truecm> point at 000 1500
  \arrow <0.15cm> [0.25,0.75] from  -1000 0000 to 4400 0000
  \arrow <0.15cm> [0.25,0.75] from 0000 -600 to  0000 4400 
  \plot 0000 0000 0400 4000 / 
  \plot 3900 50 3900 -50 / 
  \put{$a$} at 3900 -500
  \setquadratic
  \plot 0400 4000 1500 1200 4000 0200 / 
  \setdashes \setlinear
  \plot 0400 0000 0400 4000 / 
  \put{$1\over n$} at 0400 -500
  \put{$t^{\beta-\alpha}$} at 2000 1600
  \endpicture 
  }
  $$ 
  \bigskip
  For $t\geq 1/n$, it is clear that $t^\beta = g_n(t)t^\alpha$, while for
$t\leq 1/n$, we have
  $$
  \big|t^\beta - g_n(t)t^\alpha\big| =
  \big|t^\beta - n^{\alpha-\beta+1}t^{\alpha+1}\big| \ \leq \
  t^\beta + n^{\alpha-\beta+1}t^{\alpha+1} \ \leq \
  2  n^{-\beta},
  $$
  which can be made less than $\varepsilon$ for a suitably large $n$.
  \proofend

\state Lemma \label LemaOnBanachModules Let $A$ be a C*-algebra and
let $M$ be a left Banach module\fn{This means that $M$ is a Banach
space which is also a left module over $A$, satisfying
$\|am\|\leq\|a\|\|m\|$, for all $a\in A$, and $m\in M$.} over
$A$. Given a self-adjoint element $h$ in $A$ with $\|h\|\leq1$, and
$\alpha>0$, one has that
  \izitem
  \zitem 
  $
  \closure{h^\alpha M} = 
  \{m\in M: m=\lim_{k\to\infty} h^{1/k}m\},
  $
  \zitem
  If $\alpha'>0$, then $\closure{h^\alpha M} =
\closure{h^{\alpha'} M}$, 
  \zitem if $m\in \closure {h^\alpha M}$, and $h^\beta m=0$, for
some $\beta>0$, then $m=0$.

\proof We begin by proving (i). For this,  let $m\in \closure{h^\alpha M}$.
Given $\varepsilon>0$, choose $n\in M$ such that $\|m-h^\alpha n\|\leq
\varepsilon/2$.  Then, for every $k\in\N$, one has that
  \def\[{\big\Vert}\def\]{\big\Vert}%
  $$
  \[m - h^{1/k}m\] =
  \[(1 - h^{1/k})m\]\leq
  \[(1 - h^{1/k})(m-h^\alpha n)\] +   \[(1 - h^{1/k})h^\alpha
n\] \$\leq
  \[1 - h^{1/k}\]\[m-h^\alpha n\] +
  \[h^\alpha - h^{\alpha+1/k}\]\[n\] \$\leq
  {\varepsilon \over 2} + 
  \sup_{0\leq t \leq 1}   \big|t^\alpha - t^{\alpha+1/k}\big|\, \[n\],
  $$
  which can be made less than $\varepsilon$ for all sufficiently large
$k$. Observe that we have referred to ``1" above without having
assumed that $A$ is unital.  Actually our use of ``1" is just a
calculation resource which can be avoided by expanding out all
products.  This shows that
  $$
  m=\lim_{k\to\infty} h^{1/k}m,
  \subeqmark GotLimit
  $$
  as desired. 

Let $k\in\N$ be fixed. 
Plugging $\beta=1/k$ in \lcite{\ElementaryOnFunctions} we obtain a
continuous function $g$ defined on $[0,1]$ such that
  $$
  \|h^{1/k}-g(h)h^\alpha \|< \varepsilon.
  $$
  Given any $m\in M$, we then have
  $$
  \|h^{1/k}m- g(h)h^\alpha  m\| \leq
  \|h^{1/k}- g(h)h^\alpha \| \|m\| \leq \varepsilon\|m\|.
  $$
  Since $g(h)h^\alpha m = h^\alpha g(h) m$ lies in $h^\alpha M$, and $\varepsilon$ is
arbitrary, we conclude that $h^{1/k}m\in \closure{h^\alpha M}$.  If we
assume that $m$ satisfies \lcite{\GotLimit}, it then follows easily that
$m\in \closure{h^\alpha M}$, and hence (i) is proved.

Since the right-hand-side of (i) does not depend on $\alpha$, we see
that (ii) follows.

We finally prove (iii). By the same
reasoning above, given $k\in\N$ and $\varepsilon >0$, there exists a
continuous function $g$ such that
  $$
  \|h^{1/k}m\| = 
  \|h^{1/k}m-  g(h)h^\beta m\|
  \leq \varepsilon\|m\|.
  $$
  Since $\varepsilon$ is arbitrary it follows that $h^{1/k}m=0$, and
hence that $m=0$, by \lcite{\GotLimit}.
\proofend

\state Proposition \label ClosureEM 
  Let $\Cx:A\to B$ be a conditional expectation, let $M$ be a slice,
and let $\tau$ be as in \lcite{\IntroduceTau}.  Then the sets
$\tau\rg M$ and $\Cx(M)$ have the same closure.
  
\proof  \def\exponent{t}
For every $m,n\in M$ we have by \lcite{\IntroduceTauEquation} that 
  $$
  \tau mn^* = \Cx(m)\Cx(n^*) = \Cx\big(m\Cx(n^*)\big) \in \Cx(M),
  $$
  so $\tau MM^*\subseteq \Cx(M)$, and hence $\closure{\tau\rg
M}\subseteq \closure{\Cx(M)}$.
  In order to prove the reverse inclusion we first claim that 
  $$
  \Cx(u_i) \in \closure{ \tau\rg M},
  \subeqmark ANiceClaim
  $$
  for all $i$,
  where $\{u_i\}_{i\in\N}$ is any given frame for $M$.  Fixing $i\in\N$, let $\sigma
= \Cx(u_i)\Cx(u_i)^*$, and observe that a simple computation using the
C*-identity
  $$
  \big\|\sigma^{\exponent} \Cx(u_i) - \Cx(u_i)\big\|^2 = 
  \big\|\big(\sigma^{\exponent} \Cx(u_i) - \Cx(u_i)\big)
    \big(\sigma^{\exponent} \Cx(u_i) - \Cx(u_i)\big)^{\hbox{$*$}}\big\|,
  $$
  shows that
  $$
  \lim_{t\to0^+}  \sigma^{\exponent} \Cx(u_i) = \Cx(u_i).
  $$
  Clearly $\sigma\leq\tau$, so that
$\sigma^{\exponent}\leq\tau^{\exponent}$, for all $t\in(0, 1)$ by
\scite{\Pedersen}{1.3.8}, and hence
  $$ 
  \Cx(u_i)^*\sigma^{\exponent}\Cx(u_i) \leq
  \Cx(u_i)^*\tau^{\exponent}\Cx(u_i)  \leq
  \Cx(u_i)^*\Cx(u_i),
  $$
  which implies that
  $\ds
  \lim_{t\to0^+} \Cx(u_i)^*\tau^{\exponent}\Cx(u_i) = \Cx(u_i)^*\Cx(u_i).
  $
  If this time one uses the C*-identity 
  $$
  \big\|\tau^{\exponent} \Cx(u_i) - \Cx(u_i)\big\|^2 = 
  \big\|\big(\tau^{\exponent} \Cx(u_i) - \Cx(u_i)\big)^{\hbox{$*$}}
    \big(\tau^{\exponent} \Cx(u_i) - \Cx(u_i)\big)\big\|,
  $$
  there should be no difficulty in proving that $\ds\lim_{t\to0^+}
\tau^{\exponent} \Cx(u_i) = \Cx(u_i)$.  Since $\Cx(u_i)\in \rg M$, by
\lcite{\FirstLocationOfEM}, we have that 
  $$
  \tau^{\exponent} \Cx(u_i)  \in
  \tau^{\exponent} \rg M \subseteq
  \closure{\tau^{\exponent} \rg M}\  \={(\LemaOnBanachModules.ii)} \
  \closure{\tau \rg M},
  $$
  from which claim \lcite{\ANiceClaim} follows readily.  Observing that
  $\closure{\tau \rg M}$ is a {\cltsi} in $\rg M$,
  and hence also in $B$,
  and using \lcite{\GotFrame}, we then have for every $m\in M$, that
  $$
  \Cx(m) =
  \sum_{i\in \N} \Cx(u_iu_i^*m) =
  \sum_{i\in \N} \Cx(u_i)u_i^*m \in
  \closure{\tau \rg M},
  $$
  proving that $\Cx(M)\subseteq   \closure{\tau \rg M}$, and concluding
the proof. \proofend

Let us now briefly discuss the fact that a slice $M$ is a left module,
not only over $\rg M$, but also over the multiplier algebra of the
latter.  We have already mentioned that if $\{v_i\}_i$ is an
approximate identity for $\rg M$, then $\lim_iv_im = m$, for all
$m\in M$.  Given any multiplier $\mu$ of $\rg M$, and $m\in M$,
observe that
  $$
  \lim_i(\mu v_i) m 
  $$ 
  exists, since by 
  Cohen-Hewitt's factorization Theorem \scite{\CH}{32.22}
one may write $m=an$, for some $a\in \rg M$, and $n\in M$, in which
case
  $$
  \lim_i(\mu v_i) m =   \lim_i(\mu v_ia)n = (\mu a) n.
  $$
  This also shows that the correspondence 
  $$
  m=an\in M \ \longmapsto \ \mu m:= (\mu a)n\in M
  $$
  is well defined, thus making $M$ a left module over the multiplier
algebra of $\rg M$.

Regarding the element $\tau$ of the multiplier algebra of $\rg M$
produced by \lcite{\IntroduceTau}, we then have that $\tau m\in M$, 
for any $m\in M$, where the meaning of $\tau m$ is given by the above
module structure.  The same clearly also  applies to $\tau\half m$.

\state Proposition \label IntroVC
  Let $\Cx:A\to B$ be a conditional expectation, and let $M$ be a
slice.  Then there exists an isometric {\vc} $\phi$ of $B$ in $A$
whose domain is $\closure{\Cx(M)} = \closure{\tau\rg M}$, and such
that
  $$
  \phi\big(\Cx(m)\big) = \tau\half m
  \for m\in M.
  $$
  Moreover  the range of $\phi$ is precisely $\closure{\tau\half M}$.

\proof Given $m\in M$ we have 
  $$
  \|\Cx(m)\|^2 =
  \|\Cx(m)\Cx(m^*)\| \={(\IntroduceTauEquation)}
  \|\tau mm^*\| = 
  \|\tau\half mm^*\tau\half\| = 
  \|\tau\half m\|^2.
  $$
  This said one sees that there exists an isometric linear map
$\phi:\Cx(M)\to M$, such that
  $
  \phi\big(\Cx(m)\big) = \tau\half m,
  $ 
  for every $m$ in $M$, which therefore may be continuously extended to
the closure of ${\Cx(M)}$.  By abuse of language we will denoted the
extended map also by $\phi$.  It is then clear that the range of this
extended map is $\closure{\tau\half M}$.

We will now prove that $\phi$ is a {\vc}. Begining with
\lcite{\VirtuComm.i}, let $b\in B$, and $x\in \closure{\Cx(M)}$.  It
is clearly enough to consider the case in which $x=\Cx(m)$, for some
$m\in M$.  Write $m=an$, with $a\in \rg M$ and $n\in M$ as above.
Picking an approximate unit $\{v_i\}_i$ for $\rg M$, we have that
  $$
  b\phi(x) = 
  b\phi\big(\Cx(m)\big) = b (\tau\half a) n= 
  \lim_i b v_i (\tau\half a) n = 
  \lim_i \big((b v_i) \tau\half\big) a n \$= 
  \lim_i \big(\tau\half(b v_i)\big) a n = 
  \big(\tau\half(b a)\big)n =
  \phi\big(P(ban)\big) =
  \phi\big(b\Cx(m)\big) =
  \phi(bx).
  $$
  The proof of  \lcite{\VirtuComm.ii} is easier:
  $$
  \phi(x)b = \phi\big(\Cx(m)\big)b =  (\tau\half a) n b =
  \phi \big(\Cx(an b)\big) =
  \phi \big(\Cx(m)b\big) =
  \phi (xb).
  \proofend
  $$

\section{{\Gc} subalgebras}
  In this section we introduce our proposed broadening of Renault's
notion of Cartan subalgebras and prove a result on uniqueness of
conditional expectations which is a C*-algebraic analogue of
\scite{\TakeTwo}{IX.4.3} and also a generalization of
\scite{\Renault}{5.7}.

\definition \label DefineCartan Let $A$ be a C*-algebra and let
$B\subseteq A$ be a closed *-subalgebra.  We will say that $B$ is a
\stress{\gc} subalgebra of $A$ if
  \izitem
  \zitem $B$ contains an approximate unit for $A$,
  \zitem $B$ satisfies property (\MaxPrime),
  \zitem \zitemmark NormalizerInDefCartan $N(B)$ generates $A$,
  and
  \zitem there exists a
  faithful\fn{Traditionally a conditional expectation is said to be
faithful if $\Cx(a^*a)=0 \imply a=0$.  However we may adopt
a weaker form of faithfulness here by requiring only that:
$(\forall b, c\ \ \Cx(c^*ab)=0) \imply a=0$.}
  conditional expectation $\Cx: A\to B$.

Notice that $N(B)$ is closed under multiplication and under taking
adjoints, and hence its linear span is always a selfadjoint subalgebra
of $A$.  Thus, to say that
  $N(B)$ generates $A$ \stress{as a C*-algebra} is the same as saying that
  $N(B)$ generates $A$ \stress{as a Banach space}.

\state Theorem \label CondexpInCartan Let $A$ be a separable
C*-algebra, let $B\subseteq A$ be a closed *-subalgebra satisfying
\lcite{\DefineCartan.i--ii},  and let $M$ be a slice.   Given a conditional
expectation $P$ from $A$ to $B$ one has that 
  \izitem
  \zitem $M$ is invariant under $\Cx$,
  \zitem $M = \big(B\cap M\big)\oplus \big(\Ker(\Cx)\cap M\big)$,
  \zitem if $B\cap M=\0$, then $\Cx(M)=\0$.

\proof
Considering the {\vc} $\phi$ obtained from \lcite{\IntroVC} one has
by \lcite{\DefineCartan.ii}
that the range of $\phi$ is contained in $B$.  By the last sentence of
\lcite{\IntroVC} we have that $\tau\half M\subseteq B$, and hence also 
  $
  \tau M \subseteq B.
  $
  So
  $$
  \tau MM^* =
  MM^*\tau =
  M(\tau M)^* \subseteq 
  MB^* \subseteq M. 
  $$
  Using \lcite{\ClosureEM} we then deduce that
  $$
  \Cx(M) \subseteq
  \closure {\Cx(M)} =
  \closure{\tau MM^*}
  \subseteq M,
  $$
  proving (i).  Since the  restriction of $\Cx$ is an idempotent operator
on $M$ it follows that $M$ is the direct sum of $\Ker(id -\Cx|_M)$ and
$\Ker(\Cx|_M)$, from where (ii) follows.  Finally, (iii) is an immediate
consequence of (ii).  \proofend

The following is a C*-algebraic analogue of \scite{\TakeTwo}{IX.4.3}
and also a generalization of \scite{\Renault}{5.7}.

\state Theorem \label UniqueCondexp Let $A$ be a separable C*-algebra
and let $B\subseteq A$ be a closed *-subalgebra satisfying
\lcite{\DefineCartan.i--\NormalizerInDefCartan}.  Then there exists at
most one conditional expectation from $A$ onto $B$.

\proof Let $\Cx$ and $Q$ be conditional expectations onto $B$.  In
order to show that $Q=\Cx$, it is enough to show that $Q(a) = \Cx(a)$, for
every $a\in N(B)$.  Using \lcite{\SliceForEveryone} choose
a slice $M$ with $a\in M$,  and write $a=a_1+a_2$, with $a_1\in B\cap
M$, and $a_2\in \Ker(\Cx)\cap M$, according to
\lcite{\CondexpInCartan.ii}.  Since $a_1\in B$, we have that
  $$
  \Cx(a_1) = a_1 = Q(a_1),
  $$
  so it suffices to prove that $Q(a_2) =0 \ \big( =\Cx(a_2)\big)$.  

Let $N= \Ker(\Cx)\cap M$, and observe that $N$ is a slice with $B\cap
N=\0$.  Applying \lcite{\CondexpInCartan.iii} to $N$ and $Q$ we
see that $Q(N)=\0$, so $Q(a_2) =0$. \proofend

To see the relevance of $B$ satisfying property (\MaxPrime), as
opposed to just satisfying $B'\cap A\subseteq B$ (compare
\lcite{\CompareclassicalNotion}) in the above result, observe that in
the example given after \lcite{\CompareclassicalNotion} one has that
any element of $C(X)\*B$ of the form $u\*b$, where $u$ is an
unimodular function, is in $N(B)$, and hence $N(B)$ generates
$C(X)\*B$. That is, $B$ satisfies
\lcite{\DefineCartan.\NormalizerInDefCartan}. Clearly $B$ also
satisfies \lcite{\DefineCartan.i} and, although
\lcite{\DefineCartan.ii} fails, one has that $B'\cap
\big(C(X)\*B\big)\subseteq B$, as observed earlier.
  Nevertheless each probability measure on $X$ would give, by
integration, a different conditional expectation from $C(X)\*B$ to
$B$.

\section{Inverse semigroups and Fell bundles from {\gc} subalgebras}
  This section is dedicated to constructing a Fell bundle over an
inverse semigroup from an inclusion of C*-algebras.  In the case of 
{\gc} subalgebras we will later show that this Fell bundle contains
enough information to allow for a reconstruction of the inclusion.

Our standing hypothesis here will be \lcite{\StandingOne}, namely we
will assume that $A$ is a separable C*-algebra and $B$ is a closed
*-subalgebra containing an approximate unit for $A$.
This will allow us to use \lcite{\MMStarContainedInB} as in the
following result:

 \state Proposition \label DefineProdSlices Let $M$ and $N$ be slices.
Then $\overline{MN}$ (closed linear span) is a slice.

\proof
  Let $x=\sum_{i=1}^k m_in_i$, where $m_i\in M$ and $n_i\in N$.  Then,
for every $i,j$ we have   by \lcite{\MMStarContainedInB} that
  $$
  m_i^*n_i^*Bn_jm_j\subseteq M^*N^*BNM\subseteq B,
  $$
  from where we deduce that $x^*Bx\subseteq B$.  Similarly
$xBx^*\subseteq B$, so $x\in N(B)$.  This shows that
$\overline{MN}$ is contained in $N(B)$, and hence that it is a slice.
  \proofend

From now on we will frequently make use of sets such as
$\closure{MN}$, so we shall adopt the following:

\sysstate{Convention}{\rm}{\label ConventionOnProduct If $X$ and $Y$
are linear subspaces of $A$ we will denote by $X\prd Y$ the \stress{closed
linear span} of the set 
  $
  \{xy:x\in X,\ y\in Y\}.
  $}

\state Proposition \label SABIsISG The set
  $
  \SAB = \{M\subseteq A: M \hbox{ is a slice}\}
  $
  is an inverse semigroup under the operation referred to in
\lcite{\DefineProdSlices} (and from now on denoted simply by $MN$, 
according to \lcite{\ConventionOnProduct}).
  
  \proof It is evident that $\SAB$ is an associative semigroup.  Given
$M\in \SAB$ we may view $M$ as a left Hilbert module over $\rg M$
by \lcite{\SlicesTRO.i},  so we deduce that $M\prd M^*\prd M = M$ by
\scite{\Jensen}{1.1.4}, and similarly $M^*\prd M\prd
M^* = M^*$.

  Since $M^*\prd M$ and $M\prd M^*$ (now meaning closed linear span)
are {\cltsi s} in $B$ by \lcite{\SlicesTRO.ii}, they commute and hence
$\SAB$ is an inverse semigroup by \scite{\Lawson}{1.1.3}.
  \proofend

If $M$ is an idempotent element of $\SAB$, then $M= M^*\prd M$, so $M$ is a
{\cltsi} of $B$ by \lcite{\SlicesTRO.ii}.  Conversely, it is easy to see that
every {\cltsi} of $B$ is an idempotent element of
$\SAB$.  Thus
  $$
  E(\SAB) = \{J: J \hbox{ is a {\cltsi} of } B\}.
  $$
  In particular $B$ itself is an idempotent element, actually the
unit of
$\SAB$.

  The usual inverse semigroup  order on $\SAB$ is defined by 
  $$
  N\leq M \iff M\prd N^*\prd N = N.
  $$
  Notice that, when $N\leq M$, one has that
  $N = M\prd N^*\prd N \subseteq M\prd B \subseteq M$, so $N\subseteq M$.
Conversely, if $N\subseteq M$, then
  $$ 
  N =
  N\prd N^*\prd N  \subseteq
  M\prd N^*\prd N \subseteq
  M\prd M^*\prd N \subseteq
  B\prd N \subseteq N,
  $$ 
  so $N = M\prd N^*\prd N$.   In other words we have 
  $$
  N\leq M \iff N\subseteq M.
  $$

In case $N(B)$ generates $A$, as in
\lcite{\DefineCartan.\NormalizerInDefCartan}, notice that by
\lcite{\SliceForEveryone} one has that
  $$
  A = \clsum_{M\in\sSAB}M,
  \eqmark AIsSumSAB
  $$
  where we use the symbol $\clsum$ to denote the \stress{closure} of
the sum of a given family of subspaces of $A$.

\state Proposition \label GeneratingConditions
  Let $\S$ be any *-subsemigroup\fn{By a *-subsemigroup of $\sSAB$ we
shall mean a subsemigroup $\sS\subseteq\sSAB$ such that $\sS^*=\sS$, so
that $\sS$ is an inverse semigroup in its own right.} of $\SAB$.
  If $A = \clsum_{M\in\sS}M$, then $B = \clsum_{M\in\sS}M^*M$.

\proof
  Let $B' = \clsum_{M\in\sS}M^* M$, so that $B'$ is a {\cltsi} of $B$.
Also let $A'={B'AB'}$, so $A'$ is a C*-subalgebra of $A$.  For
$M\in\S$ notice that
  $$
  M = (M\prd M^*)\prd  M\prd  (M^*\prd M) \subseteq B'\prd  M\prd  B'
  \subseteq A',
  $$
  so the hypothesis implies that $A'=A$.  It is therefore clear
that any bounded approximate
unit $\{\app_i\}_i$ for $B'$ is also an approximate
unit for $A$.  Given an arbitrary $b\in B$, one therefore has that
  $$ 
  b = \lim_i\app_i b \in B'\prd B \subseteq B',
  $$
  so $B\subseteq B'$. Since the reverse inclusion also holds, the
proof is complete. \proofend

Our next result is related to one of the main hypothesis of
\scite{\actions}{9.9}. See also \scite{\actions}{5.4}.

\state Proposition \label MysteriousConditions
  Let $\S$ be any *-subsemigroup of $\SAB$.  Then the following are
equivalent:
  \izitem
  \zitem for every $M,N\in\S$, and every $a\in M\cap N$, there exists $K\in\S$
such that $a\in K\subseteq M\cap N$,
  \zitem for every $M\in\S$, and every $b\in M\cap B$, there exists $K\in\S$
such that $b\in K\subseteq M\cap B$. 

\proof (i)$\imply$(ii): Notice that $M\cap B$ is a {\cltsi} of $B$ so,
given $b\in M\cap B$, we may write $b = b_1^*b_2$, with $b_1, b_2\in
M\cap B$.  In particular $b\in M^*\prd M$.  Pluging $N = M^*\prd M$ in (i) we
deduce that there exists $K\in\S$ such that
  $$
  b\in K \subseteq M\cap M^*\prd M \subseteq M\cap B.
  $$

\noindent (ii)$\imply$(i):  Given $a\in M\cap N$, one has that
  $$ 
  a^*a \in M^*\prd M\cap M^*\prd N \subseteq B\cap M^*\prd N.
  $$ 
  By (ii) there exists $K\in\S$, such that 
  $
  a^*a \in K\subseteq B\cap M^*\prd N,
  $
  so 
  $$
  aa^*a \in MK\subseteq MB\cap M\prd M^*\prd N \subseteq M\cap N.
  $$
  The proof will therefore be concluded once we show that $a\in MK$.

  Using \lcite{\ElementaryOnFunctions}, for every $n\in\N$, let $h_n$
be a continuous real valued function on the interval $\big[0,
\|a\|^2\big]$ such that
  $$
  \big|h_n(t)t - t^{1/n}\big| < {1\over n}
  \for t\in \big[0, \|a\|^2\big],
  $$
  so that $\|h_n(aa^*)aa^* - (aa^*)^{1/n}\|<1/n$.
  We claim that 
  $$
  \lim_{n\to\infty}h_n(aa^*)aa^*a = a.
  $$
  In fact
  $$
  \|h_n(aa^*)aa^*a  - a\| \leq
  \|h_n(aa^*)aa^*a - (aa^*)^{1/n}a\| + \|(aa^*)^{1/n}a - a\| \$\leq
  {\|a\|\over n} +
    \Big\| \Big((aa^*)^{1/n}a - a\Big) \Big((aa^*)^{1/n}a -
a\Big)^{\textstyle*} \Big\|\half \$=
  {\|a\|\over n} +
    \big\| (aa^*)^{1+2/n} - 2(aa^*)^{1+1/n} + aa^*\big\|\half ,
  $$
  which converges to zero as $n\to\infty$, proving our claim.  Since
$h_n(aa^*)aa^*a \in MK$, we deduce that $a\in MK$, concluding the
proof. \proofend

\definition \label DefineAdmissible Let $\S$ be a *-subsemigroup of
$\SAB$.  We will say that $\S$ is an \stress{admissible semigroup} for
the C*-algebra inclusion ``$A\subseteq B$" if
  \izitem 
  \zitem $A = \clsum\nolimits_{M\in\sS}M$, and 
  \zitem $\S$ satisfies the equivalent conditions of
\lcite{\MysteriousConditions}.

If $N(B)$ generates $A$, we have seen that \lcite{\AIsSumSAB}
holds. Moreover $\SAB$ is clearly closed under intersections, and
hence $\SAB$ satisfies \lcite{\MysteriousConditions.i}.  In other
words, $\SAB$ is an admissible semigroup.  It is however likely to be
very big, perhaps even uncountable, so one might wonder whether a
countable admissible semigroup exists.

\state Proposition \label ExistCtblAdmiss If $(A, B)$ is a {\gc} pair
with separable $A$, then there exists a countable admissible
semigroup.

\proof If $A$ is separable then so is $N(B)$, hence we may choose a
countable dense subset $\{a_n:n\in\N\}$ of $N(B)$.  Using
\lcite{\SliceForEveryone}, for each $n\in\N$ pick a slice $M_n$
containing $a_n$.  The smallest *-subsemigroup of $\SAB$ closed under
pairwise intersections, and containing all of the $M_n$ satisfies the
conditions in the statement.  \proofend

\section{Fell bundle models for {\gc} subalgebras}
  In this section we present our main result giving a sufficient
condition for an inclusion of C*-algebras to be modeled by the
monomorphism $\Phi_\reduced: C^*(\E) \to \CstarRed(\A)$ of
\lcite{\IdentifyEinA}.

Throughout this section we let $A$ be a separable C*-algebra and $B$
be a {\gc} subalgebra of $A$.  We will also fix any
admissible semigroup $\S\subseteq\SAB$.  Observe that a countable such
$\S$ exists by \lcite{\ExistCtblAdmiss}.

We thus get a Fell bundle $\A$ over $\S$ by setting $A_M=M$, for
every $M\in\S$.  As in section \lcite{\PureStatesSection} we let $\E$
denote the restriction of $\A$ to the idempotent semilattice
$E(\S )$.  Observe that, if for each $J\in E(\S )$ we consider the
inclusion mapping
  $$
  \iota_J: J \hookrightarrow B,
  $$
  then $\iota=\{\iota_J\}_{J\in E(\sS)}$ is a representation of $\E$
in $B$ satisfying the hypothesis of \lcite{\IsomorphicIfOneToOne}.
  Since $\Ran(\iota)$ coincides with $B$ by
\lcite{\GeneratingConditions} we deduce that $C^*(\E)$ is isomorphic
to $B$.  We shell henceforth identify $C^*(\E)$ and $B$, observing
that the identification puts in correspondence the canonical mappings
  $$
  \piu_J\quad \simeq\quad \iota_J.
  \eqmark piuIsIj
  $$

Turning our attention to $\A$ let us extend the above
representation  $\iota$ to the whole of $\A$ by considering the
inclusion mappings
  $$
  \iota_M: M \hookrightarrow A,
  \eqmark HereIsRepr
  $$
  for each $M\in\A$.  It is then evident that $\iota=\{\iota_M\}_{M\in
\sS}$ is a representation of $\S$ in $A$.
  Let 
  $$
  \Phi:\LA \to A
  \eqmark HereIsMorphism
  $$
  be the *-homomorphism associated to $\iota$, as
in \lcite{\StarRepVsPreRep}.
By \lcite{\RepsVanishingOnNVsPhi} one has that $\Phi$ vanishes
on $\NA $. 

We next want to give an explicit description of the state $\tphi$ on
$\LA $ obtained from \lcite{\BigExtension} for each given pure state
$\phi$ on $C^*(\E)$.

\state Proposition \label TphiAndCondexp Let $\phi$ be a pure state on
$B = C^*(\E)$ and let $\tphi$ be the state on $\LA $ given by
\lcite{\BigExtension} in terms of $\phi$.  Then
  $$
  \tphi = \phi\circ \Cx\circ\Phi.
  $$

\proof It is clearly enough to verify that 
  $$
  \tphi(a_M\d_M) = \phi\circ \Cx\circ\Phi(a_M\d_M),
  \subeqmark ProvingStatesAgree
  $$ 
  for all $M\in \S$, and every $a_M\in M$.  

\case 1 Assume that $a_M\in\Ker(\Cx)\cap M$.  The right-hand-side 
of \lcite{\ProvingStatesAgree} then becomes
  $$
  \phi\circ \Cx\circ\Phi(a_M\d_M) = 
  \phi\big(\Cx(a_M)\big) = 0.
  $$
  Using the description of $\tphi$ given in
\lcite{\BigExtension.\CalculoDeTPhi} suppose first that that there
exists $J\in \supp(\phi)$ with $J\leq M$.  Then the left-hand-side of
\lcite{\ProvingStatesAgree} becomes
  $$
  \tphi(a_M\d_M) = \tphi_J^M(a_M) \= {(\extendingToAs)}
\lim_i\phi(a_M\app_i),
  \subeqmark ThisIsLeftHandSide
  $$
  where $\{\app_i\}_i$ is any approximate identity for $J$.  Notice
however that since $J\leq M$, we have that $MJ=J$, so
  $
  a_M\app_i \in MJ = J\subseteq B.
  $
  Thus
  $$
  a_M\app_i = \Cx(a_M\app_i) = \Cx(a_M)\app_i = 0,
  $$
  so \lcite{\ThisIsLeftHandSide} vanishes.

Suppose next that that there exists no $J\in \supp(\phi)$ with $J\leq
M$.  Then $\tphi(a_M\d_M) = 0$, by the second clause of 
\lcite{\BigExtension.\CalculoDeTPhi}, and hence
\lcite{\ProvingStatesAgree} is proved under our first case.

\case 2 Assuming that $a_M\in B\cap M$, the right-hand-side 
of \lcite{\ProvingStatesAgree} becomes
  $$
  \phi\circ \Cx\circ\Phi(a_M\d_M) = 
  \phi\big(\Cx(a_M)\big) =
  \phi(a_M).
  $$
  Recall that $\S$ is an admissible semigroup and hence that 
\lcite{\MysteriousConditions.ii} holds.  Therefore there exists
$J\in\S$ such that $a_M\in J\subseteq M\cap B$.  Since $J\subseteq B$
we see that $J$ is in $E(\S)$, and since $J\subseteq M$ we have that
$J\leq M$.  Therefore $a_M\d_M \equiv a_M\d_J$, modulo $\NA$, and
since $\tphi$ vanishes on $\NA$ by 
\lcite{\BigExtension.\TphiVanisOnN}, the 
left-hand-side of
\lcite{\ProvingStatesAgree} becomes
  $$
  \tphi(a_M\d_M) = 
  \tphi(a_M\d_J) \={(\BigExtension.ii)}
  \phi\big(\piu_J(a_M)\big) \={(\piuIsIj)}
  \phi(a_M),
  $$
  hence \lcite{\ProvingStatesAgree} is proved under case 2.

In view of \lcite{\CondexpInCartan.ii} we see that the general case of
\lcite{\ProvingStatesAgree} follows from the two cases already treated
and the proof is then complete. \proofend

We now come to our main result:

\state Theorem \label MainResult Let $A$ be a separable C*-algebra and
let $B$ be a {\gc} subalgebra of $A$.  Choose any
admissible subsemigroup $\S$ of $\SAB$, such as given by
\lcite{\ExistCtblAdmiss}, and consider the Fell bundle $\A$ over $\S$
given by $A_M=M$, for every $M\in\S$.  Then $A$ is isomorphic to
$\CstarRed(\A)$ via an isomorphism
  $
  \Psi:\CstarRed(\A) \to \A
  $
  such that 
  $$
  \Psi\circ\piur_M(m) = m
  \for M\in\S \for m\in M.
  $$
  In addition $\Psi\big(C^*(\E)\big) = B$, where $\E$ is the
restriction of $\A$ to the idempotent semilattice of $\S$.

\proof Let $\Phi:\LA \to A$ be as in
\lcite{\HereIsMorphism}.
  We first claim that 
  $$
  \|\Phi(x)\| = \tnorm x \for x\in \LA, 
  \subeqmark NormEqualityForTnormAndPhi
  $$
  where $\tnormFunction$ was defined in \lcite{\DefineReduced}.
  With this goal in mind let $\phi$ be any pure state on $B$
and consider the state $\psi = \phi\circ \Cx$ on $A$.  Also let $(\Pi_\psi,
H_\psi, \xi_\psi)$ be the GNS representation of $A$ associated to $\psi$.
Consider the map
  $$
  U: x\in \LA  \mapsto \Pi_\psi\big(\Phi(x)\big)\xi_\psi \in H_\psi.
  $$
  For $x, y\in \LA $ we have 
  $$
  \<U(x),U(y)\> = 
  \<\Pi_\psi\big(\Phi(x)\big)\xi_\psi,\Pi_\psi\big(\Phi(y)\big)\xi_\psi\> =
  \<\Pi_\psi\big(\Phi(y^*x)\big)\xi_\psi,\xi_\psi\> \$=
  \psi\big (\Phi(y^*x)\big) =
  \phi \circ \Cx\circ \Phi(y^*x) \={(\TphiAndCondexp)} \tphi(y^*x),
  $$
  where $\tphi$ is the extension of $\phi$ to $\LA $ provided by
\lcite{\BigExtension}. This says that 
  $$
  \<U(x),U(y)\> = \<x,y\>_{\tphi},
  $$
  where $\<\cdot, \cdot\>_{\tphi}$ was defined in
\lcite{\InnerProductViaphi}.  The map $U$ thus defines an isometry
(also denoted $U$ by abuse of language)
  $$
  U:H_{\tphi}\to H_\psi,
  $$
  such that $U(\hat x) = \Pi_\psi\big(\Phi(x)\big)\xi_\psi$, for every
$x\in \LA $.  By \lcite{\DefineAdmissible.i} we have that $\Phi$ has
dense image and since $\Pi_\psi$ is cyclic we deduce that $U$ is in
fact a unitary operator.  We next claim that the diagram
  $$
  \matrix{\LA  & \buildrel {\textstyle \Phi} \over \longrightarrow & A \cr\cr
  \RiLA \downarrow\quad  && \quad \downarrow \Pi_\psi\cr\cr
  B(H_{\tphi}) & \buildrel {\textstyle {\rm Ad}_U} \over \longrightarrow & B(H_\psi)}
  $$
  commutes.  Proving this amounts to checking that 
  $$
  \Pi_\psi\big(\Phi(x)\big)U|_{\hat y} = U\RiLA(x)|_{\hat y}
  \for x, y\in \LA .
  $$
  We have 
  $$
  \Pi_\psi\big(\Phi(x)\big)U|_{\hat y} = 
  \Pi_\psi\big(\Phi(x)\big)\Pi_\psi\big(\Phi(y)\big)\xi_\psi =
  \Pi_\psi\big(\Phi(xy)\big)\xi_\psi =
  U(\widehat{xy}) = U\RiLA(x)|_{\hat y}.
  $$
  This shows that our diagram indeed commutes and hence we deduce
that, for every $x\in \LA $,
  $$
  \tnorm x =
  \sup_\phi\|\RiLA(x)\| =
  \sup_\phi\|\Pi_{\psi}\big(\Phi(x)\big)\| =
  \sup_\phi\|\Pi_{\phi\circ \Cx}\big(\Phi(x)\big)\|,
  $$
  where $\phi$ runs in the set of all pure states of $B$.
  To prove \lcite{\NormEqualityForTnormAndPhi} it now suffices to
prove that 
  $$ 
  \sup_\phi\|\Pi_{\phi\circ\Cx} (a))\| =
  \|a\| \for a\in A, 
  $$
  which is in turn equivalent to proving that
$\bigoplus_\phi\Pi_{\phi\circ\Cx}$ is a faithful representation of $A$.

If $a\in A$ is nonzero then there are $b,c\in A$ such that
$\Cx(c^*ab)\neq 0$
  (see the footnote in Definition \lcite{\DefineCartan}),
  and hence $\phi\big(\Cx(c^*ab)\big)\neq0$, for some pure state
$\phi$ on $B$, which implies that $\Pi_{\phi\circ\Cx}(a)\neq0$. This
proves the faithfulness of $\bigoplus_\phi\Pi_{\phi\circ\Cx}$, and
hence concludes the proof of \lcite{\NormEqualityForTnormAndPhi}.

It is now clear the the correspondence 
  $$
  \iotaaRed(x)\in \CstarRed(\A) \longmapsto \Phi(x) \in A
  $$
  is well defined and extends to an isomorphism $\Psi: \CstarRed(\A)
\to A$, satisfying $\Psi\circ\iotaaRed=\Phi$.  That $\Psi$ is onto
follows from the  already mentioned  fact that the range of $\Phi$ is
dense in $A$.  Given $M\in\S$ we have that
  $$
  \Psi\circ \piur_M \={(\DefinePiur)}
  \Psi\circ \iotaaRed\circ \pizero_M =
  \Phi\circ \pizero_M =
  \iota_M.
  $$
  This is to say that the identity displayed in the statement holds.

Viewing $C^*(\E)$ as a subalgebra of $\CstarRed(\A)$ according to
\lcite{\IdentifyEinA}, we have that 
  $
  C^*(\E) = \clsum_{J\in E(\sS)} \piur_J(J),
  $
  so we deduce that
  $$
  \Psi\big(C^*(\E)\big)=
  \clsum_{J\in E(\sS)} \Psi\big(\piur_J(J)\big) =
  \clsum_{J\in E(\sS)} J
  \={(\GeneratingConditions)} B.
  \proofend
  $$

  \endgroup

  \begingroup
  \bigskip\bigskip 
  \font \sc = cmcsc8 \sc
  \parskip = -1pt

  Departamento de Matem\'atica 

  Universidade Federal de Santa Catarina

  88040-900 -- Florian\'opolis -- Brasil

  \eightrm r@exel.com.br

  \endgroup
  \bye